\documentclass{amsproc}

\usepackage[english]{babel}
\usepackage[utf8]{inputenc}
\usepackage{amsmath,amssymb,amsfonts}
\usepackage{url}
\usepackage[dvipsnames,usenames]{color}
\usepackage{array}
\usepackage[pdftex]{graphicx}
\usepackage{accents}
\usepackage{tabu}

\DeclareFontFamily{U}{mathx}{\hyphenchar\font45}
\DeclareFontShape{U}{mathx}{m}{n}{
      <5> <6> <7> <8> <9> <10>
      <10.95> <12> <14.4> <17.28> <20.74> <24.88>
      mathx10
      }{}
\DeclareSymbolFont{mathx}{U}{mathx}{m}{n}
\DeclareFontSubstitution{U}{mathx}{m}{n}
\DeclareMathSymbol{\bigtimes}{1}{mathx}{"91}

\newcolumntype{L}{>{\displaystyle}l}
\newcolumntype{C}{>{\displaystyle}c}
\newcolumntype{R}{>{\displaystyle}r}

\renewcommand{\tfrac}{\genfrac{}{}{}1}

\newcommand{\R}{\mathbb R}
\newcommand{\N}{\mathbb N}
\newcommand{\Z}{\mathbb Z}
\newcommand{\C}{\mathbb C}

\newcommand{\cont}{\mathcal C}
\newcommand{\D}{\mathbb D}
\renewcommand{\Im}{\mathrm{Im}}
\renewcommand{\Re}{\mathrm{Re}}
\def\E{{\mathrm{e}}}

\def\di{\partial}
\def\til{\widetilde}
\newcommand{\pois}{\mathcal P}

\def\I{\mathfrak{i}}
\newcommand{\diff}{\mathrm{d}}
\renewcommand{\bar}{\overline}
\newcommand{\och}{\quad\mathrm{and}\quad}
\newcommand{\myref}[1]{$(\ref{#1})$}

\newcommand{\sgn}{\mathrm{sgn}}

\renewcommand{\vec}[1]{\accentset{\rightharpoonup}{#1}}

\usepackage{stackengine}
\usepackage{scalerel}
\def\rlwd{.4pt}
\def\lhexbrace{\kern1pt%
\setstackgap{S}{0pt}\def\stackalignment{l}
\ThisStyle{\scalerel*{%
  \stackunder[-\rlwd]{%
    \stackon[-\rlwd]{\rule{\rlwd}{4pt}}{\rotatebox{45}{\rule{4pt}{\rlwd}}}%
  }{\rotatebox{-45}{\rule{4pt}{\rlwd}}}%
}{\SavedStyle[}}}
\def\rhexbrace{%
\setstackgap{S}{0pt}\def\stackalignment{r}
\ThisStyle{\scalerel*{%
  \stackunder[-\rlwd]{%
    \stackon[-\rlwd]{\rule{\rlwd}{4pt}}{\rotatebox{-45}{\rule{4pt}{\rlwd}}}%
  }{\rotatebox{45}{\rule{4pt}{\rlwd}}}%
}{\SavedStyle[}}\kern1pt}

\allowdisplaybreaks

\theoremstyle{definition}
\newtheorem{definicija1}{Definition}[section]
\newtheorem{primer1}[definicija1]{Example}
\newtheorem{opomba1}[definicija1]{Remark}

\theoremstyle{plain} 
\newtheorem{lema1}[definicija1]{Lemma}
\newtheorem{izrek1}[definicija1]{Theorem}
\newtheorem{trditev1}[definicija1]{Proposition}
\newtheorem{posledica1}[definicija1]{Corollary}

\numberwithin{equation}{section}

\begin{document}

\title[Herglotz-Nevanlinna functions]{An integral representation for Herglotz-Nevanlinna functions in several variables}

\author{Annemarie Luger}
\address{Annemarie Luger, Department of Mathematics, Stockholm University, 106 91 Stockholm, Sweden}
\curraddr{}
\email{luger@math.su.se}
\thanks{The authors thankfully acknowledge financial support from SSF,  grant nr.\,AM13-0011.}

\author{Mitja Nedic}
\address{Mitja Nedic, Department of Mathematics, Stockholm University, 106 91 Stockholm, Sweden}
\curraddr{}
\email{mitja@math.su.se}
\thanks{}


\date{20-9-2016}

\begin{abstract}
In this article, a characterization of the class of Herglotz-Nevan\-linna functions in $n$ variables is given in terms of an integral representation. Furthermore, alternative conditions on the measure appearing in this representation are discussed in detail. Symmetry properties induced by the integral representation are also investigated.
\end{abstract}

\maketitle

\section{Introduction}

Herglotz-Nevanlinna functions are, by definition, functions holomorphic in the poly-upper half-plane and having non-negative imaginary part. In the case of one complex variable,   this amounts  to the usual upper half-plane in $\C$,   these functions have been considered at least  a century ago. This classical case has proven itself most useful in many applications ranging from spectral theory to electromagnetic engineering. Applications of Herglotz-Nevanlinna functions in the theory of linear passive systems are in fact our main non-mathematical motivation for studying these functions in several variables. 

The great advantage of Herglotz-Nevanlinna functions in one variable is that they possess an integral representation that allows for effective treatment of these functions, as well as providing an insight into their properties. It is therefore natural to ask for an integral representation for the general case. Already around 50 years ago  Vladimirov \cite{p2vladimirov1,p2vladimirov2}  was dealing with this question. The approach in \cite{p2vladimirov1} builds upon classical distribution theory and achieves a representation of the desired form, but has the significant drawback in that it requires the assumption that the measure appearing in the integral representation is the boundary measure of the function in question and hence cannot provide a characterization of the class. 
However, the characterization that appears in \cite{p2vladimirov2} seems overly complicated when compared to the classical result in one variable or the present result, see Remark \ref{Vladimirov-characterization}.

For the case $n=2$ a characterisation with simple requirements on the measure are given in  \cite{p2luger1} . In the present paper, we further generalize the aforementioned result by showing that in fact all Herglotz-Nevanlinna functions admit an integral representation.  More precisely, in Theorem \ref{p2thm5}, we show that a holomorphic function $q$ defined in the poly-upper half-plane is a Herglotz-Nevanlinna function if and only if it can be written in the form 
\begin{equation}\label{ii}
q(\vec{z}) = a + \sum_{\ell=1}^nb_\ell z_\ell + \frac{1}{\pi^n}\int_{\R^n}K_n(\vec{z},\vec{t})\diff\mu(\vec{t}),
\end{equation}
where $a \in \R$, $\vec{b} \in [0,\infty)^n$, the kernel function $K_n$ is defined in \eqref{p2kernel} and $\mu$ is a positive Borel measure on $\R^n$ satisfying the natural growth condition \myref{p2eqthm52} and the so-called Nevanlinna condition \eqref{p2eqthm53}.

This integral representation result then allows for a detailed treatment of both the class of Herglotz-Nevanlinna functions and the class of representing measures. The main theorem does not only eliminate the drawback of the representation presented in \cite{p2vladimirov1}, but also eliminates the need for the heavy machinery of distribution theory. Our proof instead relies on a direct approach based on an integral representation for functions on the unit polydisk having non-negative real part \cite{p2vladimirov3,p2koranyi}, which is in turn built upon Cauchy's integral formula and Helly's selection principle. 

Note that the integral in \eqref{ii} is well defined even for all $z\in(\C\setminus\R)^n$. Hence the function $q$ can, by its integral representation, be extended from the poly-upper half-plane to $(\C\setminus\R)^n$. It turns out that in contrast to the one variable case these extensions might be constant even if the function $q$ in the poly-upper half-plane is not. 
The present work not only provides a generalization of our result for the case $n=2$ from \cite{p2luger1}, but also includes  an extended discussion of the Nevanlinna condition as well as consequences concerning these extensions.

The paper is organized as follows: In Section \ref{p2sec2}, we fix notations and recall some well known classical results for Herglotz-Nevanlinna functions in one variable. Section \ref{p2sec25} is devoted to a closer look of the kernel function $K_n$. The main theorem, namely the integral representation, is formulated and proved in Section \ref{p2sec3}. Conditions equivalent to condition \eqref{p2eqthm53} are given in Section \ref{p2sec4}.  In Section \ref{p2sec6}, we discuss the above named extensions with respect to symmetry and actual dependence on the variables. 

\section{Notation and prerequisites}\label{p2sec2}

The poly-upper half-plane is denoted by $\C^{+n}$ and consists of all vectors such that all entries have positive imaginary part. The unit polydisk is denoted by $\D^n$. The letter $w$ always denotes elements of $\D$, the letter $z$ denotes elements of $\C^+$, while $t$ denotes real numbers and $s$ denotes elements of the interval $[0,2\pi)$.

The focal object of our study are Herglotz-Nevanlinna functions, which are defined as follows.

\begin{definicija1}\label{p2def21}
A holomorphic function $q\colon \C^{+n} \to \C$ with non-negative imaginary part is called a \emph{Herglotz-Nevanlinna function (in $n$ variables)}.
\end{definicija1}

When $n=1$, these functions have a classical integral representation due to Nevanlinna \cite{p2nevan1}, which was presented in its current form by Cauer \cite{p2cauer}.

\begin{izrek1}[Nevanlinna]\label{p2thm21}
A function $q\colon \C^+ \to \C$ is a Herglotz-Nevanlinna function if and only if $q$ can be written as
\begin{equation}\label{p2eqthm21}
q(z) = a + b z + \frac{1}{\pi}\int_{-\infty}^\infty\left(\frac{1}{t-z} - \frac{t}{1 + t^2}\right)\diff\mu(t)
\end{equation}
where $a \in \R$, $b \geq 0$ and $\mu$ is a positive Borel measure on $\R$ satisfying
\begin{equation}\label{p2eqthm22}
\int_{-\infty}^\infty\frac{1}{1 + t^2}\diff\mu(t) < \infty.
\end{equation}
\end{izrek1}

We now mention a well known consequence of Theorem \ref{p2thm21} that will be used later on.

\begin{trditev1}\label{p2prop22}
Let $q$ be a Herglotz-Nevanlinna function in one variable. Then the following two statements hold.
\begin{itemize}
\item[(a)]{The non-tangential limit
\begin{equation}\label{p2eqprop21}
\lim\limits_{z \hat{\to} \infty}\frac{q(z)}{z} = b,
\end{equation}
where $b\geq0$ is the number that appears in representation \myref{p2eqthm21}.}
\item[(b)]{The non-tangential limit
\begin{equation}\label{p2eqprop22}
\lim\limits_{z \hat{\to} t_0}(t_0-z)q(z) = \frac{1}{\pi}\mu(\{t_0\}),
\end{equation}
where $t_0 \in \R$ and measure $\mu$ is the measure that appears in representation \myref{p2eqthm21}.}
\end{itemize}
\end{trditev1}

Recall that the notation $z \hat{\to} \infty$ denotes the limit $|z| \to \infty$ in the Stoltz domain $\{z \in \C^+~|~\theta \leq \arg(z) \leq \pi-\theta\}$ where $\theta \in (0,\tfrac{\pi}{2}]$, and similarly $z \hat{\to} t_0$ denotes the limit $z \to t_0$ in $\{z \in \C^+~|~\theta \leq \arg(z-t_0) \leq \pi-\theta\}$.

In order to effectively present and prove our main result, the introduction of certain notation is necessary. Throughout this paper, we use the convention that the accent $\vec{\color{white} a}$ denotes vectors. When this accent appears over a number, it is meant to denote a vector whose entries are all equal to the number in question. It will often be practical to use the convention that division by a vector is understood as division by the product of all the entries of that particular vector. Likewise, the multiplication by a vector is understood as multiplication by the product of all the entries of the vector in question, while the inverse of a vector is defined as by taking the inverses of the individual components of the vector. The conjugate of a vector is also defined by taking the conjugates of the entries. It is denoted by the $*$ symbol while the accent $\bar{\color{white} a}$ is still used to denote the conjugate of complex numbers. Meanwhile, a limit of a vector towards another vector is interpreted as a limit where each component of the first vector tends to the corresponding component of the second vector. Non-tangential limits of vectors are interpreted similarly.

The letters $B$ and $k$ are used throughout this paper exclusively in the following manner. The letter $B$ is used to represent subsets of the set $\{1,2,\ldots,n\}$, where $n \in \N$. For fixed $n$ and  a particular subset $B$, we denote by $B^c$ its complement and by $k$ its size. We also define $B_\ell := \{1,2,\ldots,n\}\setminus\{\ell\}$.

In connection with this notation we also introduced vectors indexed by the sets $B$ and $B^c$. Given the nature of these sets, we understand this as a grouping of the elements of the vector into two in a way specified by $B$ and $B^c$. This will often be used in conjunction with the notion of dividing by a vector to describe the situation where we divide a particular expression by all the entries of a vector that appear at positions specified by a certain set.

Finally, we introduce hexagonal brackets do denote the cases where a set $B$ does not determine the absolute ordering of the input of a function. Let us illustrate this more precisely by an example. Take $n=6$, $B=\{1,4,5\}$, $B^c=\{2,3,6\}$, $q\colon \C^{+6} \to \C$ and $\vec{z} = (z_1,z_2,z_3,z_4,z_5,z_6) \in \C^{+6}$. We then have
$$q(\I\vec{1}_B,\vec{z}_{B^c}) = q(\I,\I,\I,z_2,z_3,z_6) \quad\text{and}\quad q\lhexbrace\I\vec{1}_B,\vec{z}_{B^c}\rhexbrace = q(\I,z_2,z_3,\I,\I,z_6).$$
In the first case, the vector $\I\vec{1}_B$ comes first and the vector $\vec{z}_{B^c}$ comes second. Therefore, the first three inputs are set to $\I$ while the other three are occupied by the second, third and sixth element of $\vec{z}$. However, in the second case, even though the vector $\I\vec{1}_B$ comes first and $\vec{z}_{B^c}$ comes second, the ordering is not determined absolutely. Instead, the positions specified by the set $B$ are set to $\I$ while the positions specified by $B^c$ are set to the elements lying at the appropriate positions of the vector $\vec{z}$.

\section{The kernel $K_n$ and its properties}\label{p2sec25}

Let us begin by introducing the kernel function $K_n$, which will be a key ingredient of the integral representation formula. For $\vec{z} \in \C^{+n}$ and $\vec{t} \in \R^n$, it is defined as
\begin{equation}\label{p2kernel}
K_n(\vec{z},\vec{t}) := \I\left(\frac{2}{(2\I)^n}\prod_{\ell=1}^n\left(\frac{1}{t_\ell-z_\ell}-\frac{1}{t_\ell+\I}\right)-\frac{1}{(2\I)^n}\prod_{\ell=1}^n\left(\frac{1}{t_\ell-\I}-\frac{1}{t_\ell+\I}\right)\right).
\end{equation}
The kernel $K_n$ defined in this way is a complex constant multiple of the Schwartz kernel that appears for example in \cite{p2vladimirov1}. Note that for $n=1$, it reduces to the integrand in  the classical representation Theorem \ref{p2thm21}.

Consider now a positive Borel measure $\mu$ on $\R^n$ satisfying the growth condition
\begin{equation}\label{p2eqthm52}
\int_{\R^n}\prod_{\ell=1}^n\frac{1}{1 + t_\ell^2}\diff\mu(\vec{t}) < \infty.
\end{equation}
Given such a measure,  define the function $g\colon \C^{+n} \to \C$ as
\begin{equation}\label{p25eq7}
g\colon\vec{z} \mapsto  \frac{1}{\pi^n}\int_{\R^n}K_n(\vec{z},\vec{t})\diff\mu(\vec{t}), 
\end{equation}
which is well-defined due to \eqref{p2eqthm52} and holomorphic on $\C^{+n}$, since the kernel $K_n$ is holomorphic and locally uniformly bounded on compact subsets of $\C^{+n}$. Note that, in general, the function $g$ does not have non-negative imaginary part. For example, if $\mu$ is a point measure, then $g$ also attains values in the lower half-plane.

These functions, however, have a specific growth property that will later be used in the proof of the main theorem.

\begin{lema1}\label{p25lem1}
Let the function $g\colon\C^{+n} \to \C$ be defined by formula \myref{p25eq7} for some measure $\mu$ satisfying condition \myref{p2eqthm52}. Then
$$\lim\limits_{\vec{z}_B \hat{\to} \vec{\infty}}\frac{g(\vec{z})}{\vec{z}_B} = 0$$
for any non-empty set $B \subseteq \{1,\ldots,n\}$.
\end{lema1}

\proof
The result follows immediately from the definition of the kernel $K_n$ and Lebesgue's dominated convergence theorem.
\endproof

The imaginary part of the kernel $K_n$ will be of particular interest for us, especially because of its relation to the Poisson kernel $\pois_n$ of $\C^{+n}$. Recall that the Poisson kernel can be written using complex coordinates as
\begin{equation}\label{p2poiss}
\pois_n(\vec{z},\vec{t}) := \prod_{\ell=1}^n\frac{\Im[z_\ell]}{|t_\ell - z_\ell|^2}.
\end{equation}
Note that $\pois_n > 0$ for any $\vec{z} \in \C^{+n}$ and any $\vec{t} \in \R^n$. 

For $n=1$, it is easy to check that $\Im[K_1(z,t)] = \pois_1(z,t)$. This implies that Theorem \ref{p2thm21} also provides an integral representation of the imaginary part of a Herglotz-Nevanlinna function, which is   also used in the proof of Theorem \ref{p2thm21} to show that the imaginary part of a function defined by representation \myref{p2eqthm21} is indeed non-negative. 

For $n>1$, the corresponding equality does not hold. To determine the precise relation between $\Im[K_n]$ and $\pois_n$, we define for $j \in \{1,2,\ldots,n\}$ the expressions
$$N_{-1,j} := \frac{1}{t_j - z_j} - \frac{1}{t_j - \I},~N_{0,j} := \frac{1}{t_j - \I} - \frac{1}{t_j + \I},~N_{1,j} :=\frac{1}{t_j + \I} - \frac{1}{t_j - \bar{z_j}}.$$
We are now ready to show the following relation between $\Im[K_n]$ and $\pois_n$.

\begin{trditev1}\label{p2propK}
The equality 
\begin{equation}\label{p2eqpropK}
\Im[K_n(\vec{z},\vec{t})] = \pois_n(\vec{z},\vec{t})  -\frac{1}{(2\I)^n}\sum_{\substack{\vec{\rho} \in \{-1,0,1\}^n \\ -1\in\vec{\rho} \wedge 1\in\vec{\rho}}}N_{\rho_1,1}N_{\rho_2,2}\ldots N_{\rho_n,n}
\end{equation}
holds for any $\vec{t} \in \R^n$ and any $\vec{z} \in \C^{+n}$.
\end{trditev1}

\proof
Using the above notation, we first rewrite
\begin{multline*}
(2\I)^n\big(\Im[K_n(\vec{z},\vec{t})] - \pois_n(\vec{z},\vec{t})\big) = \prod_{j=1}^n(N_{-1,j} + N_{0,j}) + \prod_{j=1}^n(N_{0,j} + N_{1,j}) \\
- \prod_{j=1}^nN_{0,j} - \prod_{j=1}^n(N_{-1,j}+N_{0,j}+N_{1,j}).
\end{multline*}
Expanding these products yields
\begin{eqnarray*}
\prod_{j=1}^n(N_{-1,j}+N_{0,j}+N_{1,j}) & = & \sum_{\vec{\rho} \in \{-1,0,1\}^n}N_{\rho_1,1}N_{\rho_2,2}\ldots N_{\rho_n,n}, \\
\prod_{j=1}^n(N_{0,j} + N_{1,j}) & = & \sum_{\substack{\vec{\rho} \in \{-1,0,1\}^n \\ -1\not\in\vec{\rho}}}N_{\rho_1,1}N_{\rho_2,2}\ldots N_{\rho_n,n}, \\
\prod_{j=1}^n(N_{-1,j} + N_{0,j}) & = & \sum_{\substack{\vec{\rho} \in \{-1,0,1\}^n \\ 1\not\in\vec{\rho}}}N_{\rho_1,1}N_{\rho_2,2}\ldots N_{\rho_n,n}, \\
\prod_{j=1}^nN_{0,j} & = & \sum_{\substack{\vec{\rho} \in \{-1,0,1\}^n \\ -1\not\in\vec{\rho} \wedge 1\not\in\vec{\rho}}}N_{\rho_1,1}N_{\rho_2,2}\ldots N_{\rho_n,n},
\end{eqnarray*}
and hence 
\begin{multline*}
\prod_{j=1}^n(N_{-1,j} + N_{0,j}) + \prod_{j=1}^n(N_{0,j} + N_{1,j}) - \prod_{j=1}^nN_{0,j} \\
- \prod_{j=1}^n(N_{-1,j}+N_{0,j}+N_{1,j}) = -\sum_{\substack{\vec{\rho} \in \{-1,0,1\}^n \\ -1\in\vec{\rho} \wedge 1\in\vec{\rho}}}N_{\rho_1,1}N_{\rho_2,2}\ldots N_{\rho_n,n},
\end{multline*}
which finishes the proof.
\endproof

\begin{opomba1}
Observe that the sum in the right hand side of equality \myref{p2eqpropK} is empty for $n=1$, ensuring that Proposition \ref{p2propK} is valid for all $n \in \N$.
\end{opomba1}

We now consider a specific function that will appear later in the proof of the main theorem.

\begin{primer1}\label{p2ex2503}
Let $\mu = \lambda_{\R^n}$ be the Lebesgue measure on ${\R^n}$. The function $g$, defined on $\C^{+n}$ by expression \myref{p25eq7}, is then identically equal to $\I$. To see this, we first  investigate what happens when integrating $K_n$ in the $t_1$-variable with respect to the Lebesgue measure.  

We intend to use the residue theorem. Let $R > \max\{1,|z_1|\}$ and let $\Gamma_{R_-}$ be contour in $\C$ consisting of the line segment between $-R$ and $R$ and the arc $\gamma_{R_-}$, which is the lower half-circle of radius $R$ with the origin at zero. By the residue theorem, we have
$$\begin{array}{LCL}
\int_{\Gamma_{R_-}}K_n(\vec{z},\vec{t})\diff t_1 & = & -2\pi\I\lim\limits_{\xi \to -\I}K_n(\vec{z},(\xi,t_2,\dots,t_n))(\xi + \I) \\[0.35cm]
~ &= & \pi K_{n-1}(\vec{z}_{B_1},\vec{t}_{B_1}).
\end{array}$$
Recall that $B_1 = \{1,\ldots,n\}\setminus\{1\}$.

We now estimate the part of the integral over $\Gamma_{R_-}$ that runs over $\gamma_{R_-}$. Clearly,
$$\left|\int_{\gamma_{R_-}}K_n(\vec{z},(\tau,t_2,\dots,t_n))\diff\tau\right| \leq \int_{0}^{\pi}R|K_n(\vec{z},(R\E^{-\I\theta},t_2,\dots,t_n))|\diff\theta.$$
It also holds uniformly
$$\lim\limits_{R \to \infty}R|K_n(\vec{z},(R\E^{-\I\theta},t_2,\dots,t_n))| = 0.$$
In conclusion, we see that
$$\int_\R K_n(\vec{z},\vec{t})\diff t_1 = \pi K_{n-1}(\vec{z}_{B_1},\vec{t}_{B_1})$$
for any $z_1 \in \C^+$.

We may now continue by integration in the $t_2$-variable with respect to the Lebesgue measure. An analogous reasoning gives
$$\iint_{\R^2}K_n(\vec{z},\vec{t})\diff t_1 \diff t_2 = \pi^2 K_{n-2}(\vec{z}_{B},\vec{t}_{B}),$$
where the set $B = \{3,4,\ldots,n\}$. After $n$ steps in total we arrive at the constant $\pi^nK_0$. Observing that $K_0 = \I$ finishes the example.\hfill$\lozenge$
\end{primer1}

\section{The main theorem}\label{p2sec3}

The main result of this paper is the following integral representation. 

\begin{izrek1}\label{p2thm5}
A function $q\colon \C^{+n} \to \C$ is a Herglotz-Nevanlinna function if and only if $q$ can be written as
\begin{equation}\label{p2eqthm51}
q(\vec{z}) = a + \sum_{\ell=1}^nb_\ell z_\ell + \frac{1}{\pi^n}\int_{\R^n}K_n(\vec{z},\vec{t})\diff\mu(\vec{t}),
\end{equation}
where $a \in \R$, $\vec{b} \in [0,\infty)^n$, the kernel $K_n$ is defined by formula \eqref{p2kernel}, and $\mu$ is a positive Borel measure on $\R^n$ satisfying both  the growth condition \myref{p2eqthm52} and the Nevanlinna condition
\begin{equation}\label{p2eqthm53}
\sum_{\substack{\vec{\rho} \in \{-1,0,1\}^n \\ -1\in\vec{\rho} \wedge 1\in\vec{\rho}}}\int_{\R^n}N_{\rho_1,1}N_{\rho_2,2}\ldots N_{\rho_n,n}\diff\mu(\vec{t}) = 0
\end{equation}
for all $\vec{z} \in \C^{+n}$. Furthermore, the parameters $a, \vec{b}$ and $\mu$ are unique for a given function $q$.
\end{izrek1}

\begin{opomba1}
The parameters $a, \vec{b}$ and $\mu$ are sometimes collectively referred to as the \emph{data} used to represent the function $q$ and are denoted by $(a,\vec{b},\mu)$.
\end{opomba1}

\begin{opomba1}
Note that in case $n=1$, condition \myref{p2eqthm53} is empty and the classical integral representation is included as a special case in the formulation of the theorem. 
\end{opomba1}

The idea of the proof of the theorem is, as in in \cite{p2luger1}, to transform a corresponding integral representation from the polydisk to the poly-upper half-plane. However, new technical difficulties have to be tackled in the case $n>2$.

\proof
The proof is divided into several steps, starting with the uniqueness statement, which will be used later on.

{\it Step 1:} 
We will show that the parameters $a, \vec{b}$ and $\mu$ can be calculated directly from the Herglotz Nevanlinna function $q$. To start with, given a function $q$ defined by representation \eqref{p2eqthm51}, it can easily be seen that 
$$a= \Re[q(\I\vec{1})],$$
while Lemma \ref{p25lem1} implies for $\ell=1,2,\ldots,n$ that
$$b_{\ell}=\lim\limits_{y \to \infty}\frac{q\lhexbrace\I y \vec{1}_{\{\ell\}},\I \vec{1}_{B_\ell}\rhexbrace}{y}.$$
Hence, the number $a$ and the vector $\vec{b}$ are unique for a given function $q$.

Next, we show a kind of inversion formula for the measure, similar to the classical Stieltjes inversion formula in the case $n=1$. Let $\psi\colon \R^n \to \R$ be a $\cont^1$-function, such that
$$|\psi(\vec{x})| \leq C\prod_{\ell=1}^n\frac{1}{1+x_\ell^2}$$
for some constant $C \geq 0$ and all $\vec{x} \in \R^n$. We then claim that
\begin{equation}\label{p2still}
\int_{\R^n}\psi(\vec{t})\diff\mu(\vec{t})=\lim\limits_{\vec{y} \to \vec{0}^+}\int_{\R^n}\psi(\vec{x})\Im[q(\vec{x} + \I\:\vec{y})]\diff \vec{x}, 
\end{equation}
where $\mu$ is the measure from representation \myref{p2eqthm51}.

Observe first that representation \myref{p2eqthm51} implies an integral representation of $\Im[q]$, namely
$$\Im[q(\vec{z})] = \sum_{\ell=1}^nb_\ell\Im[z_\ell] + \frac{1}{\pi^n}\int_{\R^n}\Im[K_n(\vec{z},\vec{t})]\diff\mu(\vec{t}).$$
Note that we are allowed to move the imaginary part into the integral due to $\mu$ being a real measure. Proposition \ref{p2propK}, together with the fact that $\mu$ satisfies condition \eqref{p2eqthm53}, implies further
\begin{equation}\label{p2eq711}
\Im[q(\vec{z})] = \sum_{\ell=1}^nb_\ell\Im[z_\ell] + \frac{1}{\pi^n}\int_{\R^n}\pois_n(\vec{z},\vec{t})\diff\mu(\vec{t}).
\end{equation}
Using representation \myref{p2eq711}, we can rewrite the right-hand side of equality \myref{p2still} as
$$\lim\limits_{\vec{y} \to \vec{0}^+}\int_{\R^n}\psi(\vec{x})\left(\sum_{\ell=1}^n b_\ell \Im[z_\ell] + \frac{1}{\pi^n}\int_{\R^n}\pois_n(\vec{x}+\I\:\vec{y},\vec{t})\diff\mu(\vec{t})\right)\diff \vec{x}.$$
It is readily seen that the part involving the term $\sum_{\ell=1}^n b_\ell \Im[z_\ell] $ tends to $0$. What remains is the part involving the Poisson kernel, where we can use Fubini's theorem to change the order of integration. Applying Lebesgue's dominated convergence theorem allows us to change the order of the limit and the first integral. We thus arrive at
$$\frac{1}{\pi^n}\int_{\R^n}\lim\limits_{\vec{y} \to \vec{0}^+}\int_{\R^n}\psi(\vec{x})\pois_n(\vec{x}+\I\:\vec{y},\vec{t})\diff\vec{x}\:\diff\mu(\vec{t}).$$
It remains to observe that, by a well known property of the Poisson kernel, the inner integral equals 
$$\lim\limits_{\vec{y} \to \vec{0}^+}\int_{\R^n}\psi(\vec{x})\pois_n(\vec{x}+\I\:\vec{y},\vec{t})\diff\vec{x} = \pi^n\psi(\vec{t}),$$
which proves equality \myref{p2still}.
  
Suppose now that we have a Herglotz-Nevanlinna function $q$ having two integral representations of the form \myref{p2eqthm51}, one with a measure $\mu_1$ and the other with a measure $\mu_2$. Since equality \myref{p2still} holds for both $\mu_1$ and $\mu_2$, but the right-hand side of equality \myref{p2still} is the same, we obtain
$$\int_{\R^n}\psi(\vec{t})\diff\mu_1(\vec{t}) = \int_{\R^n}\psi(\vec{t})\diff\mu_2(\vec{t})$$
for all functions $\psi$ as before. But this is only possible if $\mu_1 \equiv \mu_2$. We have thus shown that the parameters appearing in representation \myref{p2eqthm51} are unique. 

{\it Step 2:} Based on Proposition \ref{p2propK}, the first direction of the equivalence is not hard to show, namely that every function defined by representation \myref{p2eqthm51} is a Herglotz-Nevanlinna function. The integral that appears in representation \myref{p2eqthm51} is a well-defined expression since the measure $\mu$ satisfies condition \myref{p2eqthm52}. It follows that the function $q$ defined in this way is holomorphic on $\C^{+n}$, since the kernel $K_n$ is holomorphic and locally uniformly bounded on compact subsets of $\C^{+n}$. 

Since the measure $\mu$ satisfies condition \myref{p2eqthm53}, it follows also that $\Im[q]$ is represented by formula \myref{p2eq711} and is thus non-negative.

{\it Step 3:} The main part of the proof is the other direction, namely to show the existence of the integral representation. We will first show formula \myref{p2eqthm51} and property \myref{p2eqthm52}. Afterwards we will focus on property \myref{p2eqthm53}. 

{\it Step 3a:} Let us begin by assuming that $q$ is a Herglotz-Nevanlinna function. We first consider the possibility that $q$ attains a real value. Then there exists a vector $\vec{z_0} \in \C^{+n}$ such that $\Im[q(\vec{z_0})] = 0$. Since $q$ is a Herglotz-Nevan\-linna function, it is holomorphic and its imaginary part $\Im[q]$ is therefore pluriharmonic. It follows now from the maximum principle for pluriharmonic functions that $\Im[q] \equiv 0$ on $\C^{+n}$. The function $q$ then has a representation of the form \myref{p2eqthm51} with $a = \Re[q(\vec{z_0})]$, $\vec{b} = \vec{0}$ and $\mu \equiv 0$. Thus, the theorem holds in this case.

{\it Step 3b:} We may now restrict ourselves to the case when $q$ does not attain a real value. As in \cite{p2luger1}, we aim to transform a corresponding result from the polydisk. To this end, we use the biholomorphism between the upper half-plane and the unit disk, defined as 
$$\varphi\colon z \mapsto \tfrac{z - \I}{z + \I} \qquad \text{with} \qquad \varphi^{-1}\colon w \mapsto \I\tfrac{1-w}{1+w}.$$
Since the function $\varphi$ is also a bijection between the sets $\R$ and $S^1\setminus\{1\}$, it facilitates the change of variables between these sets as $\tfrac{t-\I}{t+\I} = \E^{\I s}$.

Given a Herglotz-Nevanlinna function $q$, there exists a function $f$ on $\D^n$ with positive real part, such that
$$q(\vec{z}) = \I f(\varphi(z_1),\ldots,\varphi(z_n)).$$
By \cite{p2vladimirov3,p2koranyi}, any such function $f$ admits an integral representation of the form
\begin{equation}\label{p204eq9}
f(\vec{w}) = \I~\Im[f(\vec{0})] + \frac{1}{(2\pi)^n}\int_{[0,2\pi)^n}\left(2\prod_{\ell=1}^n\frac{1}{1-w_\ell\E^{-\I s_\ell}}-1\right)\diff\nu(\vec{s}),
\end{equation}
where $\nu$ is a finite positive Borel measure on $[0,2\pi)^n$ satisfying the condition that
\begin{equation}\label{p204eq5}
\int_{[0,2\pi)^n}\E^{\I m_1 s_1}\ldots\E^{\I m_n s_n}\diff\nu(\vec{s}) = 0
\end{equation}
for any multiindex $\vec{m} \in \Z^n$ with at least one positive entry and at least one negative entry. This gives an integral representation of the function $q$ of the form
\begin{equation}\label{p25eq1}
q(\vec{z}) =  -\Im[f(\vec{0})] + \frac{\I}{(2\pi)^n}\int_{[0,2\pi)^n}\left(2\prod_{\ell=1}^n\frac{1}{1-\varphi(z_\ell)\E^{-\I s_\ell}}-1\right)\diff\nu(\vec{s}).
\end{equation}
We get the first part of representation \myref{p2eqthm51} by setting $a := -\Im[f(\vec{0})] \in \R$.

Before transforming the area of integration, we note $[0,2\pi)=\{0\}\cup(0,2\pi)$ and divide accordingly the integral over $[0,2\pi)^n$ into $2^n$ pieces. In order to effectively investigate each piece separately, we use the set-notation presented in Section \ref{p2sec2}. The set $B$ is used to designate over which variables we integrate while the set $B^c$ is used to designate which variables are set to 0. We thus define
$$[0,2\pi)_B := \bigtimes_{j=1}^n\left\{\begin{array}{rcl}
\{0\} & ; & j \not\in B, \\
(0,2\pi) & ; & j \in B,
\end{array}\right. \subseteq [0,2\pi)^n.$$
Observe that the sets $[0,2\pi)_B$ and $[0,2\pi)_{B'}$ are disjoint if $B \neq B'$ and that
$$[0,2\pi)^n = \bigcup_{B \subseteq \{1,\ldots,n\}}[0,2\pi)_B,$$
where $B$ is also allowed to be empty. This allows us to rewrite representation \myref{p25eq1} as
\begin{equation}\label{p25eq1a}
q(\vec{z}) =  -\Im[f(\vec{0})] + \frac{\I}{(2\pi)^n} \sum_{\substack{B \subseteq \{1,\ldots,n\} \\ 0 \leq |B| \leq n}} \int_{[0,2\pi)_B}\left(2\prod_{\ell=1}^n\frac{1}{1-\varphi(z_\ell)\E^{-\I s_\ell}}-1\right)\diff\nu(\vec{s}).
\end{equation}

Considering first the term with $|B|=n$, that is to say, the integral over $(0,2\pi)^n$, we do the change of variables in accordance with the function $\varphi$. The area of integration thus transforms into $\R^n$ and the measure $\nu$ transforms into a measure $\mu$, related by the chosen change of variables as
$$\diff\nu(\vec{s}) = \prod_{\ell=1}^n\frac{2}{1 + t_\ell^2}\diff\mu(\vec{t}).$$
As an immediate consequence of this transformation, we see that the measure $\mu$ satisfies condition \myref{p2eqthm52} since $\nu$ is a finite measure and
$$\int_{\R^n}\prod_{\ell=1}^n\frac{1}{1 + t_\ell^2}\diff\mu(\vec{t}) = \frac{1}{2^n}\int_{(0,2\pi)^n}\diff\nu(\vec{s}) < \infty.$$
The integral in question thus becomes
$$\frac{\I}{(2\pi)^n}\int_{(0,2\pi)^n}\left(2\prod_{\ell=1}^n\frac{1}{1-\varphi(z_\ell)\E^{-\I s_\ell}}-1\right)\diff\nu(\vec{s}) = \frac{1}{\pi^n}\int_{\R^n}K_n(\vec{z},\vec{t})\diff\mu(\vec{t}),$$
where the above equality follows exclusively from symbolic manipulation of the left-hand side along with the discussed change of variables. This gives us the integral term of representation \myref{p2eqthm51}.

Integration over the corner point, corresponding to $|B| = 0$, gives
$$\frac{\I}{(2\pi)^n}\int_{\{0\}^n}\left(2\prod_{\ell=1}^n\frac{1}{1-\varphi(z_\ell)\E^{-\I s_\ell}}-1\right)\diff\nu(\vec{s}) = \frac{\I}{(2\pi)^n}\nu(\{0\}^n)\left(2\prod_{\ell=1}^n\frac{z_\ell+\I}{2\I}-1\right).$$

For $0<|B|<n$ we define the measure $\nu_B$ on $(0,2\pi)^k$, where $k = |B|$ as in Section \ref{p2sec2}, as the restriction of the measure $\nu$ to the non-zero part of $[0,2\pi)_B$, i.e.
$$\nu_B := \nu|_{[0,2\pi)_B}\lhexbrace \:\cdot\:,\vec{0}_{B_c} \rhexbrace.$$
More precisely, take for example $n=6$ and $B = \{1,4,5\}$. Considering the natural inclusion
$$\iota_B\colon (0,2\pi)^3 \hookrightarrow [0,2\pi)_B = (0,2\pi) \times \{0\} \times \{0\} \times (0,2\pi) \times (0,2\pi) \times \{0\},$$
we have that $\nu_B(U) = \nu(\iota_B(U))$ for any Borel subset $U \subseteq (0,2\pi)^3$.

The measure $\mu_B$ on $\R^k$ is then defined as the transform of the measure $\nu_B$ with respect to the change of variables induced by the function $\varphi$.
Integration over $[0,2\pi)_B$ with respect to the measure $\nu$ can then be replaced by integration over $\R^k$ with respect to the measure $\mu_B$, which also gives rise to a function of the type considered in formula \eqref{p25eq7}. More precisely, we have 
\begin{multline*}\frac{\I}{(2\pi)^n}\int_{[0,2\pi)_B}\left(2\prod_{\ell=1}^n\frac{1}{1-\varphi(z_\ell)\E^{-\I s_\ell}}-1\right)\diff\nu(\vec{s}) \\
= \frac{1}{(2\pi)^{n-k}}\left(q_B(\vec{z}_B)\cdot\prod_{\alpha \in B^c} \frac{z_\alpha+\I}{2\I}  + q_B(\I\vec{1}_B)\cdot\left(\prod_{\alpha \in B^c} \frac{z_\alpha+\I}{2\I}-1\right) \right),
\end{multline*}
where the function $q_B$ is a function of $k$ variables and is defined as
\begin{equation}\label{p25eq2}
q_B(\vec{z}_B) := \frac{1}{\pi^k}\int_{\R^k}K_k(\vec{z}_B,\vec{t}_B)\diff\mu_B(\vec{t}_B).
\end{equation}

Representation \myref{p25eq1a} thus becomes
\begin{eqnarray*}
q(\vec{z}) & =  & a + \frac{\I}{(2\pi)^n}\nu(\{0\}^n)\left(2\prod_{\ell=1}^n\frac{z_\ell+\I}{2\I}-1\right) \\
~ & ~ & + \sum_{\substack{B \subseteq \{1,\ldots,n\} \\ 1 \leq |B| \leq n-1}}\frac{1}{(2\pi)^{n-k}}\left(q_B(\vec{z}_B)\cdot\prod_{\alpha \in B^c}\frac{z_\alpha+\I}{2\I} + q_B(\I\vec{1}_B)\cdot\left(\prod_{\alpha \in B^c}\frac{z_\alpha+\I}{2\I}-1\right) \right) \\
~ & ~ & + \frac{1}{\pi^n}\int_{\R^n}K_n(\vec{z},\vec{t})\diff\mu(\vec{t}).
\end{eqnarray*}

In what follows, we are going to show that for $1 \leq |B|<n-1$, the functions $q_B$ vanish identically, whereas in case $|B|=n-1$ they have a very particular form. 

Consider now the function of one complex variable $z$ given by $q$ with all but one variable fixed. More precisely, fix a vector $\vec{\zeta} \in \C^{+n}$ and pick a position $j \in \{1,2,\ldots,n\}$. Let now $\zeta_j = z$ be interpreted as a complex variable in the upper half-plane and consider the function
$$q_j\colon z \mapsto q(\zeta_1,\ldots,\zeta_{j-1},z,\zeta_{j+1},\ldots,\zeta_n).$$
This is a Herglotz-Nevanlinna function of one variable and thus
$$\lim\limits_{z \hat{\to} \infty}\frac{q_j(z)}{z} \geq 0$$
by Proposition \ref{p2prop22}. On the other hand, by the representation above, we have
$$\begin{array}{LCL}
\multicolumn{3}{L}{\lim\limits_{z \hat{\to} \infty}\frac{q_j(z)}{z} = \lim\limits_{z \hat{\to} \infty}\frac{q(\zeta_1,\ldots,\zeta_{j-1},z,\zeta_{j+1},\ldots,\zeta_n)}{z}} \\[0.5cm]
~ &  = & \lim\limits_{z \hat{\to} \infty}\bigg(\frac{a}{z} + \frac{\I}{(2\pi)^n}\nu(\{0\}^n)\left(2\:\frac{z+\I}{2\I\:z}\prod_{\substack{\ell=1\\ \ell \neq j}}^n\frac{\zeta_\ell+\I}{2\I}-\frac{1}{z}\right) \\[0.75cm]
~ & ~ & + \frac{1}{z}\sum_{\substack{B \subseteq \{1,\ldots,n\} \\ 1 \leq |B| \leq n-1}}\frac{1}{(2\pi)^{n-k}}\left(q_B(\vec{\zeta}_B)\cdot\prod_{\alpha \in B^c}\frac{\zeta_\alpha+\I}{2\I} + q_B(\I\vec{1}_B)\cdot\left(\prod_{\alpha \in B^c}\frac{\zeta_\alpha+\I}{2\I}-1\right)\right) \\[0.75cm]
~ & ~ & + \left.\frac{1}{\pi^n}\int_{\R^n}\frac{K_n((\zeta_1,\ldots,\zeta_{j-1},z,\zeta_{j+1},\ldots,\zeta_n),\vec{t})}{z}\diff\mu(\vec{t})\right).
\end{array}$$
The limits of the first two terms are straight forward to calculate and for the last term Lemma \ref{p25lem1} implies that the limit vanishes. In the remaining sum, for all the sets $B$ such that $j \in B$, Lemma \ref{p25lem1} can be applied again and we get that
$$\lim\limits_{z \hat{\to} \infty}\left(\frac{1}{(2\pi)^{n-k}}\left(\frac{q_B(\vec{\zeta}_B)}{z}\cdot \prod_{\alpha \in B^c}\frac{\zeta_\alpha+\I}{2\I} + \frac{q_B(\I\vec{1}_B)}{z}\cdot\left(\prod_{\alpha \in B^c}\frac{\zeta_\alpha+\I}{2\I}-1\right)\right)\right) = 0.$$
Finally, for those sets $B$ for which $j \notin B$, we have
\begin{multline*}
\lim\limits_{z \hat{\to} \infty}\left(\frac{1}{(2\pi)^{n-k}}\left(\frac{z + \I}{2\I\:z}\prod_{\alpha \in B^c\setminus\{j\}}\frac{\zeta_\alpha+\I}{2\I} \cdot q_B(\vec{\zeta}_B) \right.\right. \\
\left.\left.+ \frac{1}{z}\left(\frac{z+\I}{2\I}\prod_{\alpha \in B^c\setminus\{j\}}\frac{\zeta_\alpha+\I}{2\I}-1\right)\cdot q_B(\I\vec{1}_B)\right)\right) \\
= \frac{1}{(2\pi)^{n-k}}\: \frac{1}{2\I}\big(q_B(\vec{\zeta}_B) + q_B(\I\vec{1}_B)\big)\prod_{\alpha \in B^c\setminus\{j\}}\frac{\zeta_\alpha+\I}{2\I}.
\end{multline*}
In conclusion, we arrive at
\begin{multline}\label{p2eq_b_long}
\lim\limits_{z \hat{\to} \infty}\frac{q_j(z)}{z} = \frac{\nu(\{0\}^n)}{(2\pi)^n}\prod_{\substack{\ell=1\\ \ell \neq j}}^n\frac{\zeta_\ell+\I}{2\I} \\
+ \sum_{\substack{B \subseteq \{1,\ldots,n\} \\ 1 \leq |B| \leq n-1 \\ j \notin B}}\bigg(\frac{1}{(2\pi)^{n-k}}\:\frac{1}{2\I}\big(q_B(\vec{\zeta}_B) + q_B(\I\vec{1}_B)\big) \prod_{\alpha \in B^c\setminus\{j\}}\frac{\zeta_\alpha+\I}{2\I}\bigg) = b_j
\end{multline} 
where $b_j \geq 0$. Note that in principle we have $b_j = b_j(\vec{\zeta}_{B_j})$, where $B_j = \{1,2,\ldots,n\}\setminus\{j\}$ as in Section \ref{p2sec2}. But the function $\vec{\zeta}_{B_j} \mapsto b_j(\vec{\zeta}_{B_j})$ is a holomorphic function in $n-1$ complex variables, as it is given by equality \myref{p2eq_b_long}, that is also non-negative. Therefore, it is a constant function.

Observe that, for a fixed index $j$, there is only one set $B$ in the sum above that satisfies the conditions $B \subseteq \{1,\ldots,n\}, |B| = n-1, j \notin B$. It is precisely the set $B_j$. Hence, we can express $q_{B_j}$ from equality \myref{p2eq_b_long} as
\begin{multline}\label{p2eq710}
q_{B_j}(\vec{\zeta}_{B_j}) = 4\pi\I\:b_j - q_{B_j}(\I\vec{1}_{B_j}) -  \frac{2\I\:\nu(\{0\}^n)}{(2\pi)^{n-1}}\prod_{\substack{\ell=1\\ \ell \neq j}}^n\frac{\zeta_\ell+\I}{2\I} \\
 -\sum_{\substack{B \subseteq \{1,\ldots,n\} \\ 1 \leq |B| \leq n-2 \\ j \notin B}}\bigg(\frac{1}{(2\pi)^{n-k-1}}\:\big(q_B(\vec{\zeta}_B) + q_B(\I\vec{1}_B)\big)\prod_{\alpha \in B^c\setminus\{j\}}\frac{\zeta_\alpha+\I}{2\I}\bigg).
\end{multline}
Recall that the function $q_{B_j}$ was defined by formula \eqref{p25eq2} and hence its growth behaviour is described by Lemma \ref{p25lem1}, which we are now going to use several times in order to show that the above representation can be simplified substantially. 

To start with, by  Lemma \ref{p25lem1}, it holds that
$$\lim\limits_{\vec{\zeta}_{B_j} \hat{\to} \vec{\infty}}\frac{q_{B_j}(\vec{\zeta}_{B_j})}{\vec{\zeta}_{B_j}} = 0,$$
while using expression \myref{p2eq710} of the function $q_{B_j}(\vec{\zeta}_{B_j})$ gives that
$$\lim\limits_{\vec{\zeta}_{B_j} \hat{\to} \vec{\infty}}\frac{q_{B_j}(\vec{\zeta}_{B_j})}{\vec{\zeta}_{B_j}} = \frac{2\I}{(2\pi)^{n-1}}\:\frac{1}{(2\I)^n}\nu(\{0\}^n).$$
Thus $\nu(\{0\}^n) = 0$ and the expression for $q_{B_j}$ simplifies to
\begin{multline}\label{p2eq700}
q_{B_j}(\vec{\zeta}_{B_j}) = 4\pi\I\:b_j - q_{B_j}(\I\vec{1}_{B_j}) \\
-\sum_{\substack{B \subseteq \{1,\ldots,n\} \\ 1 \leq |B| \leq n-2 \\ j \notin B}}\bigg(\frac{1}{(2\pi)^{n-k-1}}\:\big(q_B(\vec{\zeta}_B) + q_B(\I\vec{1}_B)\big)\prod_{\alpha \in B^c\setminus\{j\}}\frac{\zeta_\alpha+\I}{2\I}\bigg).
\end{multline}

We now want to show that every term in the big sum in formula \myref{p2eq700} is identically zero. 
This is done sequentially, first for all sets $B$ of size $1$, then for all sets $B$ of size $2$, and so on... and finally for all sets $B$ of size $n-2$. We describe now the process at a general step when eliminating all sets $B$ of size $d$ for some $1 \leq d \leq n-2$, where we assume that all steps corresponding to sizes $1,2,\ldots,d-1$ have already been done.

If that is the case, then formula \myref{p2eq700} has already been reduced to
\begin{multline}\label{p2eq701}
q_{B_j}(\vec{\zeta}_{B_j}) = 4\pi\I\:b_j - q_{B_j}(\I\vec{1}_{B_\eta}) \\
-\sum_{\substack{B \subseteq \{1,\ldots,n\} \\ d \leq |B| \leq n-2 \\ j \notin B}}\bigg(\frac{1}{(2\pi)^{n-k-1}}\:\big(q_B(\vec{\zeta}_B) + q_B(\I\vec{1}_B)\big)\prod_{\alpha \in B^c\setminus\{j\}}\frac{\zeta_\alpha+\I}{2\I}\bigg).
\end{multline}
Pick any set $\til{B}$ from the sum with $|\til{B}| = d$. Lemma \ref{p25lem1} says now that
$$\lim\limits_{\vec{\zeta}_{\til{B}^c\setminus\{j\}} \hat{\to} \vec{\infty}}\frac{q_{B_j}(\vec{\zeta}_{B_j})}{\vec{\zeta}_{\til{B}^c\setminus\{j\}}} = 0.$$
Calculating the same limit using formula \myref{p2eq701} instead, we get that
$$\lim\limits_{\vec{\zeta}_{\til{B}^c\setminus\{j\}} \hat{\to} \vec{\infty}}\frac{q_{B_j}(\vec{\zeta}_{B_j})}{\vec{\zeta}_{\til{B}^c\setminus\{j\}}} = \frac{1}{(2\pi)^{n-d-1}}\:\frac{1}{(2\I)^{n-d-1}}\:\big(q_{\til{B}}(\vec{\zeta}_{\til{B}}) + q_{\til{B}}(\I\vec{1}_{\til{B}})\big).$$
This follows from the fact that the only term that remains in the sum in \myref{p2eq701} is the one where $B^c\setminus\{j\}=\til{B^c}\setminus\{j\} $.
Thus
\begin{equation}\label{p2eq720}
0 = \frac{1}{(2\pi)^{n-d-1}}\:\frac{1}{(2\I)^{n-d-1}}\:\big(q_{\til{B}}(\vec{\zeta}_{\til{B}}) + q_{\til{B}}(\I\vec{1}_{\til{B}})\big).
\end{equation}
Setting first $\vec{\zeta}_{\til{B}} = \I\vec{1}_{\til{B}}$ into equality \myref{p2eq720} we get that $q_{\til{B}}(\I\vec{1}_{\til{B}}) = 0$. It follows then from equality \myref{p2eq720} that $q_{\til{B}}(\vec{\zeta}_{\til{B}}) = 0$ for any $\vec{\zeta}_{\til{B}} \in \C^{+d}$ and so $q_{\til{B}} \equiv 0$.

After the calculation of $n-2$ steps like this, we see that formula \myref{p2eq700} finally transforms into
\begin{equation}\label{p25eq3}
q_{B_j}(\vec{\zeta}_{B_j}) = 4\pi\I\:b_j - q_{B_j}(\I\vec{1}_{B_j}).
\end{equation}
By setting $\vec{\zeta}_{B_j} = \I\vec{1}_{B_j}$ into equality \myref{p25eq3}, we see that $q_{B_j}(\I\vec{1}_{B_j}) = 2\pi\I\:b_j$ and thus
\begin{equation}\label{p25eq4}
q_{B_j}(\vec{\zeta}_{B_j}) = 2\pi\I\:b_j,
\end{equation}
where $b_j$ is a non-negative number.

Summed up, since the choice of the index $j$ was arbitrary, we have obtained for $|B| = n-1$ that the corresponding function $q_B$ is constant, namely of the form \eqref{p25eq4}, while for $1 \leq |B| \leq n-2$ the functions $q_B$ are identically equal to $0$. This simplifies representation \myref{p25eq1} to the desired form, namely
$$q(\vec{z}) = a + \sum_{j = 1}^nb_j z_j+ \frac{1}{\pi^n}\int_{\R^n}K_n(\vec{z},\vec{t})\diff\mu(\vec{t}).$$

{\it Step 3c:} It is left to show that the measure $\mu$ satisfies condition \myref{p2eqthm53}. Recall that the measure $\nu$ satisfies condition \myref{p204eq5} for any multiindex $\vec{m} \in \Z^n$ such that it has at least one entry positive and at least one entry negative.

Recall also that $\vec{w}$ denotes elements in the unit polydisk $\D^n$ and define $u_j := \E^{\I s_j} \in S^1$ for $j = 1,2,\ldots,n$. Let now $\vec{\rho} \in \{-1,0,1\}^n$ and define for $j =1,2,\ldots,n$ the factors $G_{\rho_j}$ and $F_{\rho_j}$ as
$$G_{\rho_j,j} := \left\{\begin{array}{rcl}
w_j\bar{u_j}& ; & \rho_j = -1, \\[0.25cm]
1 & ; & \rho_j = 0, \\[0.25cm]
\bar{w_j}u_j & ; & \rho_j = 1,
\end{array}\right. \och F_{\rho_j,j} := \left\{\begin{array}{rcl}
\frac{w_j\bar{u_j}}{1-w_j\bar{u_j}} & ; & \rho_j = -1, \\[0.25cm]
1 & ; & \rho_j = 0, \\[0.25cm]
\frac{\bar{w_j}u_j}{1-\bar{w_j}u_j} & ; & \rho_j = 1,
\end{array}\right.$$
respectively.

If $\vec{\rho} \in \{-1,0,1\}^n$ is such that $1 \in \vec{\rho}$ and $-1 \in \vec{\rho}$ (i.e. both $1$ and $-1$ appear at least once in the vector $\vec{\rho}$), then by condition \myref{p204eq5} we have 
$$\int_{[0,2\pi)^n}G_{\rho_1,1}^{m_1}G_{\rho_2,2}^{m_2}\ldots G_{\rho_n,n}^{m_n}\diff\nu(\vec{s}) = 0$$
and hence also 
$$\int_{[0,2\pi)^n}\sum_{\vec{m} \in \N^n}G_{\rho_1,1}^{m_1}G_{\rho_2,2}^{m_2}\ldots G_{\rho_n,n}^{m_n}\diff\nu(\vec{s}) = 0,$$
since the geometric series permits the interchange of the order of integration and summation. Note that due to the condition on $\vec{\rho}$ the sum in the above expression is indeed well defined (in particular, not all of the factors $G_{\rho_j,j}$ are equal $1$).

Since
$$\sum_{\vec{m} \in \N^n}G_{\rho_1,1}^{m_1}G_{\rho_2,2}^{m_2}\ldots G_{\rho_n,n}^{m_n} = F_{\rho_1,1}F_{\rho_2,2}\ldots F_{\rho_n,n},$$
it holds that
$$\int_{[0,2\pi)^n} F_{\rho_1,1}F_{\rho_2,2}\ldots F_{\rho_n,n}\diff\nu(\vec{s}) = 0$$
and, in particular,
\begin{equation}\label{p2eq600}
\sum_{\substack{\vec{\rho} \in \{-1,0,1\}^n \\ -1 \in \vec{\rho} \wedge 1 \in \vec{\rho}}}\int_{[0,2\pi)^n}  F_{\rho_1,1}F_{\rho_2,2}\ldots F_{\rho_n,n}\diff\nu(\vec{s}) = 0.
\end{equation}

In what follows, we show that the area of integration in formula \myref{p2eq600} can be reduced to the open cube $(0,2\pi)^n$. To do that, pick any index $\vec{\rho}$ in the sum. For that particular index we can write the integral over $[0,2\pi)^n$ as
$$\int_{[0,2\pi)^n}  F_{\rho_1,1}F_{\rho_2,2}\ldots F_{\rho_n,n}\diff\nu(\vec{s}) = \sum_{B \subseteq \{1,\ldots,n\}}\int_{[0,2\pi)_B}  F_{\rho_1,1}F_{\rho_2,2}\ldots F_{\rho_n,n}\diff\nu(\vec{s}).$$

When $1 \leq |B|<n-1$, we know already from the previous step that the corresponding function $q_B$ was identically equal to zero. Hence, by the uni\-queness result of step {\it Step 1} of the proof, we have that the measure $\mu_B$ is the zero measure. This implies in turn that the measures $\nu_B$ are also zero measures and hence the terms in the sum corresponding to sets $B$ of this particular size vanish. This holds even in the case $B=\emptyset$, as we have shown in \textit{Steb 3b} that $\nu(\{0\}^n) = 0$. 

For $|B|=n-1$, the corresponding function $q_B$ was a real multiple of $\I$. In Example \ref{p2ex2503}, we have seen that the Lebesgue measure is a representing measure for the function $\I$. Therefore, again by the uniqueness result of {\it Step 1}, the measure $\mu_B$ is a constant multiple of the Lebesgue measure $\lambda_{\R^{n-1}}$ and the corresponding term in the above sum can be calculated explicitly. To start with, applying the change of variables given by the function $\varphi$, we find  that
$$F_{\rho_j,j} \xrightarrow{\varphi} \frac{1+t_j^2}{2\I}N_{\rho_j,j}.$$
Hence, for the Lebesgue measure, the relevant term in the above sum factorizes, where at least one factor corresponds to an index $\varrho_j\not=0$. It is straight forward to check that 
$$\int_{\mathbb R} N_{-1,j}\, \diff \lambda_{\R}(t_j)= \int_{\R} \left(\frac1{t_j-z_j}-\frac1{t_j-\I}\right) \diff t_j = 0$$
for $z_j \in \C^+$. Similarly,
$$\int_{\mathbb R} N_{1,j}\, \diff \lambda_{\R}(t_j)= \int_{\R} \left(\frac1{t_j+\I}-\frac1{t_j-\bar{z_j}}\right) \diff t_j = 0$$
for $z_j \in \C^+$. This implies that, indeed, in the above sum, only the term that corresponds to integration over $(0,2\pi)^n$ remains, as all other terms were shown to be individually equal to zero.  We therefore conclude that
\begin{multline*}
0=\sum_{\substack{\vec{\rho} \in \{-1,0,1\}^n \\ -1 \in \vec{\rho} \wedge 1 \in \vec{\rho}}}\int_{(0,2\pi)^n}  F_{\rho_1,1}F_{\rho_2,2}\ldots F_{\rho_n,n}\diff\nu(\vec{s}) \\
 = \frac{1}{\I^n}\sum_{\substack{\vec{\rho} \in \{-1,0,1\}^n \\ -1 \in \vec{\rho} \wedge 1 \in \vec{\rho}}}\int_{\R^n} N_{\rho_1,1}N_{\rho_2,2}\ldots N_{\rho_n,n}\diff\mu(\vec{t}).
 \end{multline*}
This shows that the measure $\mu$ satisfies condition \myref{p2eqthm53}, finishing the proof.
\endproof

\begin{opomba1}\label{p2rem2}
When $n=2$, the sum in equality \myref{p2eq700} is empty and the sequential process used to show that all the terms in the sum are zero is already completed before it even starts. Otherwise, the sum is non-empty for $n \geq 3$ and we thus require $n-2$ steps to show that all terms in the sum are in fact equal to zero.
\end{opomba1}

\begin{opomba1}\label{Vladimirov-characterization}
As mentioned in the introduction, also in \cite{p2vladimirov2}  a characterization of Herglotz-Nevanlinna functions via an integral representations   is obtained by transforming the corresponding result from the polydisk. However, there it is not used that the measure on the boundary of  $[0,2\pi)^n$ is either zero or a multiple of the Lebesgue-measure. Hence both in the integral representation and the condition on the measure all the lower dimensional integrals are still present. 
\end{opomba1}

In the following corollary, we want to highlight some facts which are direct consequences of the proof of Theorem \ref{p2thm5}.
\begin{posledica1}\label{p2posfour}
Let $n \geq 2$ and let $q$ be a Herglotz-Nevanlinna function. Then the following statements hold.
\begin{itemize}
\item[(i)]{If there exists a vector $\vec{z_0} \in \C^{+n}$ such that $\Im[q(\vec{z_0})] = 0$, then $q(\vec z) = q(\vec{z_0})$ for all  $\vec z\in\C^{+n}$.}
\item[(ii)]{The imaginary part of $q$ can be represented as
\begin{equation}\label{p25eq8}
\Im[q(\vec{z})] = \sum_{\ell=1}^nb_\ell\Im[z_\ell] + \frac{1}{\pi^n}\int_{\R^n}\pois_n(\vec{z},\vec{t})\diff\mu(\vec{t}),
\end{equation}
where $\vec{b}$ and $\mu$ are as in Theorem \ref{p2thm5} and $\pois_n$ is the Poisson kernel defined by \eqref{p2poiss}.}
\item[(iii)]{The number $a$ from Theorem \ref{p2thm5} is equal to
$$a = \Re[q(\I\vec{1})].$$
}
\item[(iv)]{For every set $B$ with $|B| \geq 1$, it holds that
$$\lim\limits_{\vec{z}_B \hat{\to} \vec{\infty}}\frac{q(\vec{z})}{\vec{z}_B} = \left\{\begin{array}{rcl}
b_j & ; & B = \{j\}, \\
0 & ; & \mathrm{else},
\end{array}\right.
$$
where the numbers $b_j$ are as in Theorem \ref{p2thm5}. In particular, the limit is independent of $\vec{z}_{B^c}$.}
\item[(v)]{For every set $B$ with $|B| \geq 1$, it holds that
$$\lim\limits_{\vec{y}_B \to \vec{\infty}}\frac{\Im[q(\I\:\vec{y})]}{\vec{y}_B} = \left\{\begin{array}{rcl}
b_j & ; & B = \{j\}, \\
0 & ; & \mathrm{else},
\end{array}\right.
$$
where the numbers $b_j$ are as in Theorem \ref{p2thm5}. In particular, the limit is independent of $\vec{y}_{B^c}$.}
\item[(vi)]{There exist non-positive numbers $c_j$ such that for every set $B$ with $|B| \geq 1$, it holds that
$$\lim\limits_{\vec{z}_B \hat{\to} \vec{0}}\vec{z}_B\:q(\vec{z}) = \left\{\begin{array}{rcl}
c_j & ; & B = \{j\}, \\
0 & ; & \mathrm{else}.
\end{array}\right.
$$
In particular, the limit is independent of $\vec{z}_{B^c}$.}
\item[(vii)]{For every set $B$ with $|B| \geq 1$, it holds that
$$\lim\limits_{\vec{y}_B \to \vec{0}^+}\vec{y}_B\:\Im[q\lhexbrace\I\:\vec{y}_B,\vec{z}_{B^c}\rhexbrace] = \left\{\begin{array}{rcl}
-c_j & ; & B = \{j\}, \\
0 & ; & \mathrm{else},
\end{array}\right.
$$
where the numbers $c_j$ are the same ones that appear in statement (vi). In particular, the limit is independent of $\vec{z}_{B^c}$.}
\item[(viii)]{
Let $\psi\colon \R^n \to \R$ be a $\cont^1$-function, such that
$$|\psi(\vec{x})| \leq C\prod_{\ell=1}^n\frac{1}{1+x_\ell^2}$$
for some constant $C \geq 0$ and all $\vec{x} \in \R^n$. Then
\begin{equation*}\label{p2still2}
\int_{\R^n}\psi(\vec{t})\diff\mu(\vec{t})=\lim\limits_{\vec{y} \to \vec{0}^+}\int_{\R^n}\psi(\vec{x})\Im[q(\vec{x} + \I\:\vec{y})]\diff \vec{x}.
\end{equation*}
}
\end{itemize}
\end{posledica1}

\proof
Statements (i),(ii),(iii),(iv) and (viii) have already been proven during the course of the proof of Theorem \ref{p2thm5}. Statement (v) follows immediately from statement (ii). Statement (vi) follows by considering the Herglotz-Nevanlinan function $\vec{z} \mapsto q(\vec{z}_{B^c},-\vec{z}_B^{-1})$ and then using statement (iv).

In order to show statement (vii), note that, by statement (vi), we have
$$\lim\limits_{\vec{y}_B \to \vec{0}^+}\I\:\vec{y}_B\:q\lhexbrace\I\:\vec{y}_B,\vec{z}_{B^c}\rhexbrace = x$$
for some $x \in \R$. Multiplying this equality by $-\I$ and then taking the imaginary finishes the proof.
\endproof

We finish this section with an example of a Herglotz-Nevanlinna function in three variables.

\begin{primer1}
Let
$$q(z_1,z_2,z_3) = 1 - \frac{1}{z_1+z_2+z_3}.$$
The function $q$ is then indeed a Herglotz-Nevanlinna. The parameters $a$ and $\vec{b}$ corresponding to $q$ in the sense of representation \myref{p2eqthm51} can easily be calculated using Corollary \ref{p2posfour} and are equal to
$$a = 1, \quad \vec{b} = \vec{0}.$$
The identity 
$$\Im[q(z_1,z_2,z_3)] = \frac{\Im[z_1]+\Im[z_2]+\Im[z_3]}{|z_1+z_2+z_3|^2},$$
suggests that the measure $\mu$ equals $\pi$ times the two-dimensional Lebesgue measure supported on the hyperplane $ t_1+t_2+t_3=0 $. Indeed, this  can be verified by a direct calculation.\hfill$\lozenge$
\end{primer1}

\section{Alternative descriptions  of the class of representing measures}\label{p2sec4}

In this section, we focus on condition \myref{p2eqthm53} and its alternative descriptions and mention also two further properties of the class of representing measures.

The main result of this section is the following theorem that establishes three alternative descriptions of the Nevanlinna condition.

\begin{izrek1}\label{p2thmalt}
Let $n \geq 2$ and let $\mu$ be a positive Borel measure on $\R^n$ satisfying the growth condition \myref{p2eqthm52}. Then the following statements are equivalent.
\begin{itemize}
\item[(a)]{
$$\sum_{\substack{\vec{\rho} \in \{-1,0,1\}^n \\ -1\in\vec{\rho} \wedge 1\in\vec{\rho}}}\int_{\R^n}N_{\rho_1,1}N_{\rho_2,2}\ldots N_{\rho_n,n}\diff\mu(\vec{t}) = 0$$
for all $\vec{z} \in \C^{+n}$.}\vspace{1mm}
\item[(b)]{
$$\int_{\R ^n}N_{\rho_1,1}N_{\rho_2,2}\ldots N_{\rho_n,n}\diff\mu(\vec{t}) = 0$$
for all $\vec{z} \in \C^{+}$ and for every vector $\vec{\rho} \in \{-1,0,1\}^n$ with at least one entry equal to $1$ and at least one entry equal to $-1$.}\vspace{1mm}
\item[(c)]{
$$\int_{\R ^n}\left(\frac{t_1-\I}{t_1+\I}\right)^{m_1}\ldots\left(\frac{t_n-\I}{t_n+\I}\right)^{m_n}\prod_{\ell=1}^n\frac{1}{1+t_\ell^2}\diff\mu(\vec{t}) = 0$$
for all multiindices $\vec{m} \in \Z^n$ with at least one positive entry and at least one negative entry.}\vspace{1mm}
\item[(d)]{
$$\int_{\R^n}\frac{1}{(t_{j_1}-z_{j_1})^2(t_{j_2}-\bar{z_{j_2}})^2}\:\prod_{\substack{\ell=1 \\ \ell \neq j_1,j_2}}^n\left(\frac{1}{t_\ell-z_\ell} - \frac{1}{t_\ell-\bar{z_\ell}}\right)\diff\mu(\vec{t}) = 0$$
for all $\vec{z} \in \C^{+n}$ and for all indices $j_1,j_2 \in \N$, such that $1 \leq j_1 < j_2 \leq n$.}

\end{itemize}
\end{izrek1}

\proof
The implication (a) $\Rightarrow$ (b) is contained within \textit{Step 3c} of the proof of Theorem \ref{p2thm5}. The converse implication is trivial. 

The transformation of condition \myref{p204eq5} to $\R^{n}$ using the change of variables fostered by the function $\varphi$, along with \textit{Step 3b} and \emph{Step 3c} of the proof of Theorem \ref{p2thm5}, implies also the equivalence between (a) and (c). 

The equivalence between (a) and (d) is given in   Proposition \ref{p2plurihar5} below.
\endproof

In order to complete the proof of Theorem \ref{p2thmalt} we consider first integral representations for the imaginary part of  Herglotz-Nevanlinna functions. 

\begin{trditev1}\label{p2plurihar4}
Let $n \geq 2$ and let $\mu$ be a positive Borel measure on $\R^n$ satisfying the growth condition \myref{p2eqthm52}. Then there exists a Herglotz-Nevanlinna function $q$ with
\begin{equation}\label{p26eq1}
\Im[q(\vec{z})] = \frac{1}{\pi^n}\int_{\R^n}\pois_n(\vec{z},\vec{t})\diff\mu(\vec{t})
\end{equation}
if and only if for all indices $j_1,j_2 \in \N$, such that $1 \leq j_1 < j_2 \leq n$, we have
\begin{equation}\label{p26eq2}
\int_{\R^n}\frac{1}{(t_{j_1}-z_{j_1})^2(t_{j_2}-\bar{z_{j_2}})^2}\:\prod_{\substack{\ell=1 \\ \ell \neq j_1,j_2}}^n\left(\frac{1}{t_\ell-z_\ell} - \frac{1}{t_\ell-\bar{z_\ell}}\right)\diff\mu(\vec{t}) = 0
\end{equation}
for any  vector $\vec{z} \in \C^{+n}$.
\end{trditev1}

\proof
Suppose first that there exists a Herglotz-Nevanlinna function $q$ such that its imaginary part satisfies formula \myref{p26eq1}. Since $q$ is then holomorphic, we have that $\Im[q]$ is pluriharmonic. Thus, for any two indices $j_1,j_2$ such that $1 \leq j_1,j_2 \leq n$, we have that
$$\frac{\di^2}{\di z_{j_1} \di \bar{z_{j_2}}}\int_{\R^n}\pois_n(\vec{z},\vec{t})\diff\mu(\vec{t}) \equiv 0$$
when $\vec{z} \in \C^{+n}$. Since this holds for arbitrary indexes $j_1,j_2$ with $1 \leq j_1,j_2 \leq n$, it holds also when $1 \leq j_1 < j_2 \leq n$. It is easily seen that the Poisson kernel allows for differentiation to be moved under the integral sign. Thus 
$$0 \equiv \frac{\di^2}{\di z_{j_1} \di \bar{z_{j_2}}}\int_{\R^n}\pois_n(\vec{z},\vec{t})\diff\mu(\vec{t}) = \int_{\R^n}\frac{\di^2}{\di z_{j_1} \di \bar{z_{j_2}}}\pois_n(\vec{z},\vec{t})\diff\mu(\vec{t}).$$
But the rightmost integral is nothing but a constant multiple of expression \myref{p26eq2}.

Conversely, suppose the measure $\mu$ is such that condition \myref{p26eq2} is satisfied. Then, we can define a function $v$ on $\C^{+n}$ as
$$v(\vec{z}) = \frac{1}{\pi^n}\int_{\R^n}\pois_n(\vec{z},\vec{t})\diff\mu(\vec{t}).$$
The function $v$ is pluriharmonic since $\mu$ satisfies condition \myref{p26eq2}, thus there exists a holomorphic function $q$ on $\C^{+n}$ such that $\Im[q] = v$. But since $\pois_n \geq 0$ for any $\vec{z} \in \C^{+n}$ and $\vec{t} \in \R^n$ and $\mu$ is a positive measure, then it follows that $q$ is a Herglotz-Nevanlinna function.
\endproof

\begin{trditev1}\label{p2plurihar5}
Let $n \geq 2$ and let $\mu$ be a positive Borel measure on $\R^n$ satisfying the growth condition \myref{p2eqthm52}. The measure $\mu$ then satisfies condition \myref{p2eqthm53} if and only if it satisfies condition \myref{p26eq2}.
\end{trditev1}

\proof
Assume first that the measure $\mu$ satisfies condition \myref{p2eqthm53}. By Theorem \ref{p2thm5}, there exists a Herglotz-Nevanlinna function $q$ having representation \myref{p2eqthm51} with $a = 0$, $\vec{b} = \vec{0}$ and this particular measure $\mu$. The imaginary part of $q$ then has representation \myref{p25eq8}. Proposition \ref{p2plurihar4} now implies that the measure $\mu$ satisfies condition \myref{p26eq2}.

Conversely, assume that the measure $\mu$ satisfies condition \myref{p26eq2}. By Proposition \ref{p2plurihar4}, there exists a Herglotz-Nevanlinna function $q$ such that equality \myref{p26eq1} holds. But, by Theorem \ref{p2thm5}, the function $q$ also admits a representation of the form \myref{p2eqthm51} with some parameters $a$, $\vec{b}$ and $\sigma$ as described by the theorem.

Observe now that, for $B = \{j\}$, the limit
$$\lim\limits_{\vec{y}_B \to \vec{\infty}}\frac{\Im[q(\I\:\vec{y})]}{\vec{y}_B}$$
is equal to $b_j$ by Corollary \ref{p2posfour}(v) and equal to zero by equality \myref{p26eq1} and the fact that
$$\lim\limits_{\vec{y}_B \to \vec{\infty}}\frac{1}{\pi^n}\int_{\R^n}\frac{\pois_n(\vec{x} + \I\:\vec{y},\vec{t})}{\vec{y}_B}\diff\mu(\vec{t}) = 0$$
for every set $B$ with $|B| \geq 1$. Thus $\vec{b} = \vec{0}$.

We know now that
$$\frac{1}{\pi^n}\int_{\R^n}\pois_n(\vec{z},\vec{t})\diff\mu(\vec{t}) =\frac{1}{\pi^n}\int_{\R^n}\pois_n(\vec{z},\vec{t})\diff\sigma(\vec{t})$$
and claim that $\mu = \sigma$. To show this, we need to repeat a few of the steps of the proof of uniqueness of the representing measure from Theorem \ref{p2thm5}. Note that we cannot use the uniqueness result, since we do not know that the function $q$ admits a representation of the from \myref{p2eqthm51} with the measure $\mu$ (we only know that $q$ admits a representation of the from \myref{p2eqthm51} with the measure $\sigma$).

Take now a function $\psi$ as in Corollary \ref{p2posfour}(viii). Considering the expression
$$\lim\limits_{\vec{y} \to \vec{0}^+}\int_{\R^n}\psi(\vec{x})\Im[q(\vec{x} + \I\:\vec{y})]\diff \vec{x}$$
with our two representations of $\Im[q]$ given by formulas \myref{p25eq8} and \myref{p26eq1}, we arrive at
\begin{multline*}
\frac{1}{\pi^n}\int_{\R^n}\lim\limits_{\vec{y} \to \vec{0}^+}\int_{\R^n}\psi(\vec{x})\pois_n(\vec{x} + \I\:\vec{y},\vec{t})\diff \vec{x}\:\diff\mu(\vec{t}) \\
= \frac{1}{\pi^n}\int_{\R^n}\lim\limits_{\vec{y} \to \vec{0}^+}\int_{\R^n}\psi(\vec{x})\pois_n(\vec{x} + \I\:\vec{y},\vec{t})\diff \vec{x}\:\diff\sigma(\vec{t}).
\end{multline*}
Here, we first used Fubini's theorem to switch the order of integration, followed by an application of Lebesgue's dominated convergence theorem in order to interchange the limits and the first integrals. Using a well known property of the Poisson kernel, we can simplify the limits and inner integrals in the last equality to arrive at
$$\int_{\R^n}\psi(\vec{t})\diff\mu(\vec{t}) = \int_{\R^n}\psi(\vec{t})\diff\sigma(\vec{t}).$$
Since $\psi$ was an arbitrary function form Corollary \ref{p2posfour}(viii), it follows that $\mu = \sigma$. In particular, this implies that the measure $\mu$ satisfies condition \myref{p2eqthm51} since the measure $\sigma$ satisfies condition \myref{p2eqthm51}.
\endproof

This  establishes now all equivalences described by Theorem \ref{p2thmalt} and we finish this section by presenting two properties of the class of representing measures.

\begin{trditev1}\label{p2posMeasure}
Let $n \geq 2$ and let $\mu$ be the representing measure of a Herglotz-Nevanlinna function. The following statements hold.
\begin{itemize}
\item[(i)]{The measure $\mu$ cannot be finite unless 0.}
\item[(ii)]{All point have zero mass.}
\end{itemize}
\end{trditev1}

We omit the proof and note that it is completely analogous to the case $n=2$ presented in \cite{p2luger1}. Other properties of the class of representing measures will be presented as an independent topic in an upcoming work.

\section{Symmetries}\label{p2sec6}

Theorem \ref{p2thm21} implies a symmetry property of Herglotz-Nevanlinna functions in one variable. Any such function, extended to $\C\setminus\R$ using representation \myref{p2eqthm21}, satisfies the equality
$$\bar{q(z)} = q(\bar{z})$$
for any $z \in \C\setminus\R$. Thus, if we know the values of the function in one half-plane, we automatically know them also in the other half-plane.

In what follows, we investigate the symmetries possessed by Herglotz-Nevanlinna functions in several variables as induced by Theorem \ref{p2thm5}. More precisely, we show a fascinating explicit formula that relates the values of any function defined using formula \eqref{p25eq7} between the different connected components of $(\C\setminus\R)^n$. In the case when we are working with a Herglotz-Nevanlinna function, we can even say something about the variable-dependence of the values of the function in different connected components of $(\C\setminus\R)^n$.

We begin by considering a function $g$ of the form \myref{p25eq7}. Recall that such a function is a priori defined only on the poly-upper half-plane $\C^{+n}$. However, it can be extended to $(\C\setminus\R)^n$, more precisely, under the assumption that the measure $\mu$ satisfies the growth condition \eqref{p2eqthm52}, the expression
$$\int_{\R^n}K_n(\vec{z},\vec{t})\diff\mu(\vec{t})$$
is well-defined for any $\vec{z} \in (\C\setminus\R)^n$. We illustrate this with the following example.

\begin{primer1}\label{p2ex2502}
Let $\mu = \lambda_{\R^n}$ the Lebesgue measure on ${\R^n}$. We have shown in Example \ref{p2ex2503} that the function $g$, defined on $\C^{+n}$ by expression \myref{p25eq7}, is identically equal to $\I$. Extending this function to $(\C\setminus\R)^n$ gives the function
\begin{equation}\label{p2eqfunI}
g(\vec{z}) = \left\{\begin{array}{rcl}
\I & ; & \vec{z} \in \C^{+n}, \\
-\I & ; & \textrm{else}.
\end{array}\right.
\end{equation}
To see this, we again begin by using the residue theorem to investigate what happens when we integrate $K_n$ in the $t_1$-variable with respect to the Lebesgue measure. 

If $z_1 \in \C^+$, then the calculations performed in Example \ref{p2ex2503} apply and give
$$\int_\R K_n(\vec{z},\vec{t})\diff t_1 = \pi K_{n-1}(\vec{z}_{B_1},\vec{t}_{B_1}).$$
Suppose now that $z_1 \in \C^-$. Similarly, by the residue theorem, we have
$$\begin{array}{LCL}
\int_{\Gamma_{R_+}}K_n(\vec{z},\vec{t})\diff t_1 & = & 2\pi\I\lim\limits_{\xi \to \I}K_n(\vec{z},(\xi,t_2,\dots,t_n))(\xi - \I) \\[0.35cm]
~ &= & -\pi K_{n-1}(\I\vec{1}_{B_1},\vec{t}_{B_1}).
\end{array}$$
Here, $\Gamma_{R_+}$ is the contour in $\C$ consisting of the line segment between $-R$ and $R$ and the arc $\gamma_{R_+}$, which is the upper half-circle of radius $R$ with the origin at zero. 

An analogous estimation argument to the one done in Example \ref{p2ex2503} allows us to conclude that
$$\int_\R K_n(\vec{z},\vec{t})\diff t_1 = -\pi K_{n-1}(\I\vec{1}_{B_1},\vec{t}_{B_1}).$$
Therefore, we have
$$\int_{\R}K_n(\vec{z},\vec{t})\diff t_1 = \left\{\begin{array}{rcl}
\pi K_{n-1}(\vec{z}_{B_1},\vec{t}_{B_1}) & ; & z_1 \in \C^{+}, \\[0.25cm]
-\pi K_{n-1}(\I\vec{1}_{B_1},\vec{t}_{B_1}) & ; & z_1 \in \C^{-}.
\end{array}\right.$$

We may now continue by integration in the $t_2$-variable with respect to the Lebesgue measure. An analogous reasoning gives
$$\iint_{\R^2}K_n(\vec{z},\vec{t})\diff t_1 \diff t_2 = \left\{\begin{array}{rcl}
\pi^2 K_{n-2}(\vec{z}_{B_{1,2}},\vec{t}_{B_{1,2}}) & ; & (z_1,z_2) \in \C^{+2}, \\[0.25cm]
-\pi^2 K_{n-2}(\I\vec{1}_{B_{1,2}},\vec{t}_{B_{1,2}}) & ; & \text{else},
\end{array}\right.$$
where the set $B_{1,2} = \{3,4,\ldots,n\}$. After $n$ steps in total we arrive at $\pi^nK_0$ if all of the numbers $z_\ell$ lie in $\C^+$, otherwise we arrive at $-\pi^nK_0$. Noting that $K_0 = \I$ gives formula \myref{p2eqfunI}.\hfill$\lozenge$
\end{primer1}

We continue by identifying the symmetries possessed by a function $g$ of the form \myref{p25eq7}. First, we introduce the following notation. Take any vector $\vec{z} \in (\C \setminus\R)^n$. We associate four sets of indices with this particular vector, namely $C_+$, $I_+$, $I_-$ and $C_-$. There are pairwise disjoint subsets of $\{1,\ldots,n\}$, defined by
$$\begin{array}{rclcrcl}
\ell \in C_+ & \iff & z_\ell \in \C^+\setminus\{\I\}, & \qquad & \ell \in I_+ & \iff & z_\ell = \I, \\
\ell \in C_- & \iff & z_\ell \in \C^-\setminus\{-\I\}, & \qquad & \ell \in I_- & \iff & z_\ell = -\I.
\end{array}$$
Note also that $C_+ \cup I_+ \cup I_- \cup C_- = \{1,\ldots,n\}$. Next, for a set $B$, we introduce a map $\Psi_B\colon \C^n \to \C^n$, defined as $\psi_B(\vec{z}) := \vec{\zeta}$, where
$$\zeta_\ell = \left\{\begin{array}{rcl}
\bar{z_\ell} & ; & \ell \in B, \\
\I & ; & \text{else}.
\end{array}\right.$$

Let us illustrate this notation with the following example. Take $n=6$ and 
$$\vec{z} = (\I,-\I,1+2\I,\I,-5-\I,2-\I) \in (\C\setminus\R)^6.$$
Then $C_+ = \{3\}$, $I_+ = \{1,4\}$, $I_- = \{2\}$ and $C_- = \{5,6\}$. We also have, for example,
$$\Psi_{\{1,2,5\}}(\vec{z}) = (-\I,\I,\I,\I,-5+\I,\I) \quad\text{and}\quad \Psi_{\{3,4\}}(\vec{z}) = (\I,\I,1-2\I,-\I,\I,\I).$$

We will now show the following proposition.

\begin{trditev1}\label{p26propS}
Let $n \geq 2$ and let a function $g\colon(\C\setminus\R)^n \to \C$ be defined using formula \myref{p25eq7} for some measure $\mu$ satisfying condition \myref{p2eqthm52}. Then, for any point $\vec{z} \in (\C\setminus\R)^n$, we have
\begin{equation}\label{p25eq9}
g(\vec{z}) = (\sgn|I_-|-1)\sum_{D\subseteq C_-\cup C_+}(-1)^{|D|}\left(\bar{g(\Psi_{D}(\vec{z}))}+\bar{g(\I\vec{1})}\right)+\bar{g(\I\vec{1})},
\end{equation}
where the sets $C_+,C_-,I_+,I_-$ are associated to the point $\vec{z}$ as before.
\end{trditev1}

\proof
Let $\vec{z} \in (\C\setminus\R)^n$ be arbitrary and introduce the notation
$$f(z,t) := \frac{(z+\I)(t-\I)}{2\I(t-z)},$$
where $z \in \C\setminus\R$ and $t \in \R$. Observe also that 
\begin{equation}\label{p2eq:symm}
\bar{f(z,t)} + f(\bar{z},t) = 1
\end{equation}
for any such numbers $z$ and $t$. Using this notation, we can rewrite the kernel $K_n$ with regards to the sets  $C_+,C_-,I_+,I_-$, which are determined by our starting vector $\vec{z}$. We get
\begin{multline*}
K_n(\vec{z},\vec{t}) = \I\prod_{\ell =1}^n\frac{1}{1+t_\ell^2}\left(2\prod_{\ell=1}^n\frac{(z_\ell+\I)(t_\ell-\I)}{2\I(t_\ell-z_\ell)}-1\right) \\[0.25cm]
= \I\prod_{\ell = 1}^n\frac{1}{1+t_\ell^2}\left(2\prod_{\ell \in C_+}f(z_\ell,t_\ell) \cdot \prod_{\ell \in C_-}f(z_\ell,t_\ell) \cdot \prod_{\ell \in I_+}f(z_\ell,t_\ell) \cdot \prod_{\ell \in I_-}f(z_\ell,t_\ell) -1 \right).
\end{multline*}

If $I_- \neq \emptyset$, then
$$\prod_{\ell \in I_-}f(z_\ell,t_\ell) = 0$$
since $f(-\I,t) = 0$ for any $t \in \R$. Thus, in this case, we have that
$$g(\vec{z}) = -\frac{\I}{\pi^n}\int_{\R^n}\prod_{\ell = 1}^n\frac{1}{1+t_\ell^2}\diff\mu(\vec{t}) = -g(\I\vec{1}) = \bar{g(\I\vec{1})}$$
and the proposition holds in this case.

Therefore, we may assume for the remainder of the proof that $I_- = \emptyset$. Observe now that $f(\I,t) = 1$ for any $t \in \R$. Using formula \eqref{p2eq:symm}, we may now rewrite the kernel $K_n$ as
\begin{eqnarray*}
K_n(\vec{z},\vec{t}) & = &  \I\prod_{\ell = 1}^n\frac{1}{1+t_\ell^2}\left(2\prod_{\ell \in C_+}f(z_\ell,t_\ell) \cdot \prod_{\ell \in C_-}f(z_\ell,t_\ell) -1 \right) \\
~ & = & \bar{-\I\prod_{\ell = 1}^n\frac{1}{1+t_\ell^2}\left(2\prod_{\ell \in C_+}\bar{f(z_\ell,t_\ell)} \cdot \prod_{\ell \in C_-}\bar{f(z_\ell,t_\ell)} -1 \right)} \\
~ & = & \bar{-\I\prod_{\ell = 1}^n\frac{1}{1+t_\ell^2}\left(2\prod_{\ell \in C_+}(1-f(\bar{z_\ell},t_\ell)) \cdot \prod_{\ell \in C_-}(1-f(\bar{z_\ell},t_\ell)) -1 \right)} \\
~ & = & \bar{-\I\prod_{\ell = 1}^n\frac{1}{1+t_\ell^2}\left(2\sum_{D \subseteq C_+ \cup C_-}(-1)^{|D|}\prod_{\ell \in D}f(\bar{z_\ell},t_\ell) -1 \right)} \\
~ & = & \bar{-\sum_{D \subseteq C_+ \cup C_-}(-1)^{|D|}\I\prod_{\ell = 1}^n\frac{1}{1+t_\ell^2}2\prod_{\ell \in D}f(\bar{z_\ell},t_\ell) + \I\prod_{\ell = 1}^n\frac{1}{1+t_\ell^2}} \\
~ & = & \bar{-\sum_{D \subseteq C_+ \cup C_-}(-1)^{|D|}\I\prod_{\ell = 1}^n\frac{1}{1+t_\ell^2}\left(\left(2\prod_{\ell \in D}f(\bar{z_\ell},t_\ell) -1\right) + 1\right)+\I\prod_{\ell = 1}^n\frac{1}{1+t_\ell^2}}.
\end{eqnarray*} 
Using the last expression for the kernel $K_n$ in formula \eqref{p25eq7} finishes the proof.
\endproof

We illustrate this proposition with one example.

\begin{primer1}\label{p26ex1}
Let $\mu$ be a positive Borel measure on $\R^2$ defined as
$$\mu(U) := \int_{\R^2}\frac{\chi_U(t_1,t_2)}{(1+t_1^2)(1+t_2^2)}\diff t_1 \diff t_2.$$
Observe that this measure satisfies condition \myref{p2eqthm52}. However, it does not satisfy condition \myref{p2eqthm53}. Define now a function $g$ of the form \myref{p25eq7} using this particular measure $\mu$. Explicitly, using the residue theorem, we have
$$g(z_1,z_2) = \frac{1}{\pi^2}\int_{\R^2}\frac{K_2\big((z_1,z_2),(t_1,t_2)\big)}{(1+t_1^2)(1+t_2^2)}\diff t_1 \diff t_2 = -\frac{7\I +z_1+z_2 + \I\:z_1z_2}{8(z_1+\I)(z_2+\I)}.$$
In order to consider the extension of the function to $(\C\setminus\R)^2$, we, for convenience, write variables lying in $\C^-$ as conjugates of variables lying in $\C^+$.

The extension of the function $g$ to $(\C\setminus\R)^2$ can likewise be calculated using the residue theorem and gives
\begin{eqnarray*}
g(\bar{z_1},z_2) & = &  -\frac{5\I +\bar{z_1}+3z_2 + \I\:\bar{z_1}z_2}{8(\bar{z_1}-\I)(z_2+\I)}, \\
g(z_1,\bar{z_2}) & = & -\frac{5\I +3z_1+\bar{z_2} + \I\:z_1\bar{z_2}}{8(z_1+\I)(\bar{z_2}-\I)}, \\
g(\bar{z_1},\bar{z_2}) & = & -\frac{-\I +3\bar{z_1}+3\bar{z_2} + \I\:\bar{z_1}\bar{z_2}}{8(\bar{z_1}-\I)(\bar{z_2}-\I)}.
\end{eqnarray*}
Suppose now that we have a point $\vec{z} \in (\C\setminus\R)^2$, such that $C_+ = \{2\},C_-=\{1\},I_+=I_-=\emptyset$. By formula \myref{p25eq9}, we then have
$$g(\bar{z_1},z_2) = \bar{g(z_1,\I)} + \bar{g(\I,\bar{z_2})} - \bar{g(z_1,\bar{z_2})}.$$
For our function $g$, this can even be verified explicitly, as
$$\begin{array}{LCL}
\multicolumn{3}{L}{ \bar{g(z_1,\I)} + \bar{g(\I,\bar{z_2})} - \bar{g(z_1,\bar{z_2})}} \\[0.25cm]
~ & = & \bar{\left(-\frac{1}{2(z_1+\I)}\right)} + \bar{\left(-\frac{1}{2(\bar{z_2}-\I)}\right)} - \bar{\left(-\frac{5\I +3z_1+\bar{z_2} + \I\:z_1\bar{z_2}}{8(z_1+\I)(\bar{z_2}-\I)}\right)} \\[0.5cm]
~ & = & \frac{-1}{2(\bar{z_1}-\I)} + \frac{-1}{2(z_2+\I)} + \frac{-5\I + 3\bar{z_1}+z_2-\I\bar{z_1}z_2}{8(\bar{z_1}-\I)(z_2+\I)} \\[0.5cm]
~ & = & -\frac{5\I +\bar{z_1}+3z_2 + \I\:\bar{z_1}z_2}{8(\bar{z_1}-\I)(z_2+\I)} = g(\bar{z_1},z_2).
\end{array}$$
Points $\vec{z} \in (\C\setminus\R)^2$ that have different sets $C_+$, $C_-$, $I_+$ and $I_-$ can be treated similarly.\hfill$\lozenge$
\end{primer1}

Proposition \ref{p26propS} shows that, unlike in the one-variable case, the values of a Herglotz-Nevanlinna function in $(\C\setminus\R)^n$ are not completely determined only by the values in $\C^{+n}$ unless we are interested in the value of the function at a point that has at least one coordinate equal to $-\I$. Otherwise, we require knowledge of the values of the function in $2^n -1$ connected components of $(\C\setminus\R)^n$ in order to determine the values of the function in the remaining connected component of  $(\C\setminus\R)^n$. However, the following improvement is possible.

Let us now consider a function $q_0$ of the form \eqref{p25eq7} given by a measure $\mu$ that satisfies the growth condition \myref{p2eqthm52} and that also satisfies the Nevanlinna condition \myref{p2eqthm53}. In therms of Theorem \ref{p2thm5}, we have a Herglotz-Nevanlinna function represented by $(0,\vec{0},\mu)$. We will now show that the Nevanlinna condition improves the symmetry formula \myref{p25eq9} in the following way.

\begin{trditev1}\label{p2posR}
Let $n \geq 2$ and let $q_0$ be a Herglotz-Nevanlinna function represented by $(0,\vec{0},\mu)$ in the sense of Theorem \ref{p2thm5}. If we have some $z_j \in \C^-$, then the value $q_0(\vec{z})$ does not depend on the variables which lie in $\C^+$.
\end{trditev1}

\proof
Choose any $j_+ \in C_+$ and any $j_- \in C_-$. Using Lebesgue's dominated convergence theorem, we first calculate that
$$\frac{\di}{\di z_{j_+}}q_0(\vec{z}) = \frac{1}{\pi^n}\frac{\bar{z_{j_-}}+\I}{(2\I)^{n-1}} \int_{\R^n} \frac{1}{(t_{j_+}-z_{j_+})^2(t_{j_-}-\bar{z_{j_-}})(t_{j_-} + \I)}\cdot F \cdot \diff\mu(\vec{t}),$$
where the expression $F$ is equal to
$$F = \prod_{\ell \in C_+\setminus\{j_+\}}\left(\frac{1}{t_\ell - z_\ell} - \frac{1}{t_\ell + \I}\right) \prod_{\ell \in C_-\setminus\{j_-\}}\left(\frac{1}{t_\ell - \bar{z_\ell}} - \frac{1}{t_\ell + \I}\right).$$
Observe also that we use the same convention for writing variables that lie in $\C^-$ as in Example \ref{p26ex1}.

We know from Corollary \myref{p2posfour}(viii) that $\mu$ is the distributional boundary value of $\Im[q_0]$ and can therefore write
\begin{multline*}
\frac{\di}{\di z_{j_+}}q_0(\vec{z}) \\
 =\lim\limits_{\vec{y} \to \vec{0}^+}\frac{1}{\pi^n}\frac{\bar{z_{j_-}}+\I}{(2\I)^{n-1}} \int_{\R^n} \frac{1}{(t_{j_+}-z_ {j_+})^2(t_{j_-}-\bar{z_{j_-}})(t_{j_-} + \I)}\cdot F \cdot \Im[q_0(\vec{t} + \I\:\vec{y})]\diff\vec{t}.
\end{multline*}
Writing out the imaginary part, we further get that
\begin{multline*}
\frac{\di}{\di z_{j_+}}q_0(\vec{z}) \\
 =\lim\limits_{\vec{y} \to \vec{0}^+}\frac{\bar{z_{j_-}}+\I}{(2\pi\I)^{n}} \int_{\R^n} \frac{1}{(t_{j_+}-z_{j_+})^2(t_{j_-}-\bar{z_{j_-}})(t_{j_-} + \I)}\cdot F \cdot q_0(\vec{t} + \I\:\vec{y})\diff\vec{t} \\
 - \lim\limits_{\vec{y} \to \vec{0}^+}\bar{\frac{z_{j_-}-\I}{(-2\pi\I)^{n}} \int_{\R^n} \frac{1}{(t_{j_+}-\bar{z_{j_+}})^2(t_{j_-}-z_{j_-})(t_{j_-} - \I)}\cdot \bar{F} \cdot q_0(\vec{t} + \I\:\vec{y})\diff\vec{t}}.
\end{multline*}
Denote now by $\tau$ a complexified $t$-variable and observe that the integrand in first of the above integrals has no singularities in the upper half-plane in the $\tau_{j_-}$-variable, while the integrand in the second integral has no singularities in the upper half-plane in the $\tau_{j_+}$-variable. Thus, by the residue theorem, both integral are equal to zero and so 
$$\frac{\di}{\di z_{j_+}}q_0(\vec{z}) \equiv 0.$$

Note that this argument is valid due to the fact that the function $\vec{\tau} \mapsto q_0(\vec{\tau} + \I\:\vec{y})$ is holomorphic in the domain $\bigtimes_{\ell = 1}^n\{\tau_\ell \in \C~|~\Im[\tau_\ell] > -y_\ell\}$ and the observation that
$$\lim\limits_{R \to \infty}\int_{0}^{\pi}\frac{q_0\lhexbrace\vec{\tau}_{B_{j_-}} + \I\:\vec{y}_{B_{j_-}},(R\E^{\I\theta}+\I\:y_{j_-})\vec{1}_{\{j_-\}}\rhexbrace}{(R\E^{\I\theta}-\bar{z_{j_-}})(R\E^{\I\theta}+\I)}\cdot R\:\I\:\E^{\I\theta}\diff\theta = 0$$
due to Corollary \myref{p2posfour}(iv) and the fact that the integrand in the above integral can be bounded from above independently of $\theta \in [0, \pi]$. A similar argument holds when integrating in the $\tau_{j_+}$-variable in the second integral.
\endproof

Any Herglotz-Nevanlinna function $q$ thus satisfies at any point $\vec{z} \in (\C\setminus\R)^n$ the symmetry formula
\begin{multline}\label{p25eq91}
q(\vec{z}) = a + \sum_{\ell \in C_+}b_\ell z_\ell + \sum_{\ell \in I_+}\I\:b_\ell- \sum_{\ell \in I_-}\I\:b_\ell + \sum_{\ell \in C_-}b_\ell \bar{z_\ell}  \\
+ (\sgn|I_-|-1)\sum_{D\subseteq C_-\cup C_+}(-1)^{|D|}\left(\bar{q_0(\Psi_{D}(\vec{z}))}+\bar{q_0(\I\vec{1})}\right)+\bar{q_0(\I\vec{1})},
\end{multline}
where the sets $C_+$, $C_-$, $I_+$ and $I_-$ are associated to the point $\vec{z}$. Note also that if the function $q$ is represented by $(a,\vec{b},\mu)$ in the sense of Theorem \ref{p2thm5}, then the function $q_0$ is represented by $(0,\vec{0},\mu)$ in the same sense.

We now illustrate these results with one example.

\begin{primer1}\label{p26ex2}
Let
$$q(z_1,z_2) = 2z_2 - \frac{1}{z_1+z_2}.$$
Then $q$ is a Herglotz-Nevanlinna function in two variables, where the corresponding function $q_0$ is given by
$$q_0(z_1,z_2) = -\frac{1}{z_1+z_2}.$$
The extension of the function $q_0$ to the set $(\C\setminus\R)^2$ can be calculated using the residue theorem and gives
$$\begin{array}{RCL}
q_0(\bar{z_1},z_2) & = & \frac{1}{\I-\bar{z_1}}, \\[0.25cm]
q_0(z_1,\bar{z_2}) & = & \frac{1}{\I-\bar{z_2}}, \\[0.25cm]
q_0(\bar{z_1},\bar{z_2}) & = & \frac{1}{\I-\bar{z_1}} + \frac{1}{\I-\bar{z_2}} + \frac{1}{\bar{z_1}+\bar{z_2}}.
\end{array}$$
For a point $\vec{z} \in (\C\setminus\R)^2$ with $C_+ = \{1\},C_-=\{2\},I_+=I_-=\emptyset$, formula \eqref{p25eq91} gives
$$q(z_1,\bar{z_2}) = 2\bar{z_2} + \bar{q_0(\I,z_2)} + \bar{q_0(\bar{z_1},\I)} - \bar{q_0(\bar{z_1},z_2)}.$$
For our function $q$, this can also be verified explicitly, as
$$\begin{array}{LCL}
\multicolumn{3}{L}{ 2\bar{z_2} + \bar{q_0(\I,z_2)} + \bar{q_0(\bar{z_1},\I)} - \bar{q_0(\bar{z_1},z_2)}} \\[0.25cm]
~ & = & 2\bar{z_2} + \bar{\left(-\frac{1}{\I+z_2}\right)} +  \bar{\left(\frac{1}{\I-\bar{z_1}}\right)} -  \bar{\left(\frac{1}{\I-\bar{z_1}}\right)} \\[0.25cm]
~ & = & 2 \bar{z_2} - \frac{1}{-\I+\bar{z_2}} = 2\bar{z_2} + q_0(z_1,\bar{z_2}) = q(z_1,\bar{z_2}).
\end{array}$$
Points $\vec{z} \in (\C\setminus\R)^2$ that have different sets $C_+$, $C_-$, $I_+$ and $I_-$ can be treated similarly.\hfill$\lozenge$
\end{primer1}

\section*{Acknowledgements}

The authors would like to thank Ragnar Sigurdsson for interesting discussions on the subject as well as careful reading of the manuscript.

\bibliographystyle{amsplain}

\end{document}